# The fundamental weighted category of a weighted space (*)

# (From directed to weighted algebraic topology)

Marco Grandis

*Dipartimento di Matematica, Università di Genova, Via Dodecaneso 35, 16146-Genova, Italy.*
*e-mail:* `grandis@dima.unige.it`
*home page:* `http://www.dima.unige.it/~grandis/`

**Abstract.** We want to investigate 'spaces' where paths have a 'weight', or 'cost', expressing length, duration, price, energy, etc. The weight function is *not* assumed to be invariant up to path-reversion. Thus, 'weighted algebraic topology' can be developed as an enriched version of directed algebraic topology, where illicit paths are penalised with an infinite cost, and the licit ones are measured. Its algebraic counterpart will be 'weighted algebraic structures', equipped with a sort of directed seminorm.

In the *fundamental weighted category* of a generalised metric space, introduced here, each homotopy class of paths has a weight (or seminorm), which is subadditive with respect to composition. We also study a more general setting, *spaces with weighted paths*, which has finer quotients and strong links with noncommutative geometry. Weighted homology of weighted cubical sets has already been developed in a previous work, with similar results.

## Introduction

The recent domain of directed algebraic topology studies 'directed spaces' (preordered topological spaces, locally preordered spaces, cubical sets, etc.) with 'directed algebraic structures' produced by homotopy or homology functors: on the one hand the fundamental *category* (with its higher dimensional versions), on the other *preordered* homology group. Its general aim is modelling non-reversible phenomena. See [G3, G4, G7] and references therein.

We want to propose an enrichment of this subject, *weighted algebraic topology*, replacing the truth-valued approach of directed algebraic topology (where a path is licit or not) with a measure of costs, taking values in the interval $[0, \infty]$ of positive real extended number. The general aim is, now, measuring the cost of (possibly non-reversible) phenomena.

---

(*) Work supported by MIUR Research Projects.

Weighted algebraic topology will study 'weighted spaces', like (generalised) metric spaces, with 'weighted' algebraic structures, like the fundamental *weighted* (or normed) category, defined here, and the *weighted* homology groups already developed in [G6] for weighted (or normed) cubical sets.

Lawvere's generalised metric spaces [Lw], endowed with a possibly non-symmetric distance taking values in $[0, \infty]$, are a basic setting where weighted algebraic topology can be developed (see Section 1). The present theory is essentially based on the *standard generalised metric interval* $\delta\mathbf{I}$, with distance $\delta(x, y) = y - x$, if $x \leq y$, and $\delta(x, y) = \infty$ otherwise; and the resulting cylinder functor $I(X) = X \otimes \delta\mathbf{I}$, where the tensor product has the $l_1$-type metric (Section 2). We define the fundamental weighted category $w\Pi_1(X)$ of a generalised metric spaces, and begin its study (Sections 3, 4).

As with various situations in homotopy theory, we have to work with *elementary* and *extended* homotopies, produced by 1-Lipschitz maps (i.e., weak contractions) and Lipschitz maps, respectively (see Section 2). The first are used to define the main homotopical constructs, namely cylinder, cone, suspension and - dually - cocylinder, cocone, loop-object (in the pointed case); and then to obtain the (co)fibration sequence of a map. But extended homotopies can be concatenated, and are essential to reach the higher order properties of such sequences. Also, in the fundamental weighted category itself, an arrow is a class of *extended* paths, up to *extended* homotopy with fixed endpoints.

(The reader can think of the following analogy: when studying homotopies of chain algebras, the multiplicative homotopies are used to define the main constructs, including the fibration sequence of a morphism; but such homotopies cannot be concatenated, and we need the 'extended' homotopies of the underlying chain complexes to prove the homotopy equivalence properties of the sequence, cf. [G2].)

We also introduce, in Sections 5-6, the more flexible setting of *w-spaces*, or *spaces with weighted paths*, which - like weighted cubical sets in [G6] - is able to express topological facts usually investigated with noncommutative geometry and missed by ordinary topology.

We only treat, in Section 7, an example, linked with the well-known *irrational rotation C*-algebras* $A_\vartheta$ ($\vartheta$ irrational) classified by Pimsner - Voiculescu [PV] and Rieffel [Ri], also known as 'noncommutative tori'. The algebra $A_\vartheta$ is 'meant' to replace the topological quotient $\mathbf{R}/G_\vartheta$ of the euclidean line modulo the action of the dense additive group $G_\vartheta = \mathbf{Z} + \vartheta\mathbf{Z}$, which is trivial (i.e., has the coarse topology). Here, the $\delta$-metric analogue $(\delta\mathbf{R})/G_\vartheta$ is trivial as well, with the null distance. But we get an interesting quotient from the *standard w-line* $w\mathbf{R}$, which assigns a finite weight $w(a) = a(1) - a(0)$ to each increasing path $a: \mathbf{I} \to \mathbf{R}$, and $w(a) = \infty$ otherwise (5.4). Now, the w-space $W_\vartheta = (w\mathbf{R})/G_\vartheta$ has a non-trivial fundamental weighted monoid (at any point), isomorphic to the additive monoid $G_\vartheta^+ = G_\vartheta \cap \mathbf{R}^+$ with the natural weight $w(x) = x$. We prove in Thm. 7.3, 7.4 that the *irrational rotation w-space* $W_\vartheta$ has the same classification up to isometric isomorphism (resp. Lipschitz isomorphism) as the C*-algebra $A_\vartheta$ up to isomorphism (resp. strong Morita equivalence).

**1. Generalised metric spaces**

We begin with recalling the structure of Lawvere's generalised metric spaces [Lw], which will be used as a first setting for weighted algebraic topology. Among the new points, notice the standard spaces of 1.5, the *reflective* symmetric distance (1.6) and the associated *symmetric topology* (1.9).

**1.1. Real weights.** The basic ingredient is the strict symmetric monoidal closed category of extended positive real numbers, introduced by Lawvere [Lw], which we write $\mathbf{w}^+$. It has objects $\lambda \in [0, \infty]$, morphisms $\lambda \geq \mu$, and tensor product $\lambda + \mu$ (with $\lambda + \infty = \infty$, for all $\lambda$).

As a complete lattice, this category has all limits and colimits, 'reduced' to products and sums

(1) $\quad$ product: $\sup(\lambda_i) = \vee \lambda_i,$ $\qquad$ terminal object: 0,

$\qquad$ sum: $\quad \inf(\lambda_i) = \wedge \lambda_i,$ $\qquad$ initial object: $\infty$.

The internal hom is given by *truncated subtraction*, written as a difference:

(2) $\quad \lambda + \mu \geq \nu \iff \lambda \geq \hom^+(\mu, \nu) = \nu - \mu,$ $\qquad\qquad (\nu - \mu = 0 \vee (\nu - \mu)).$

Let $\mathbf{v}$ denote the full subcategory of $\mathbf{w}^+$ on the objects $0, \infty$; in this subcategory, the cartesian product $\lambda \vee \mu$ coincides with the tensor product $\lambda + \mu$. Thus, the covariant embedding of the boolean algebra $\mathbf{2} = (\{0, 1\}, \leq)$ of *truth-values* (contravariant with respect to the natural orders)

(3) $\quad$ M: $\mathbf{2} \to \mathbf{w}^+,$ $\qquad$ M(0) $= \infty$, M(1) $= 0$,

is strict monoidal with respect to the cartesian product and the additive tensor product as well. Moreover, M has left and right adjoint

(4) $\quad$ P $\dashv$ M $\dashv$ Q, $\qquad$ P($\lambda$) = 1 $\iff \lambda < \infty,$ $\qquad\qquad$ Q($\lambda$) = 1 $\iff \lambda = 0$.

A function w: $A \to [0, \infty]$ defined on a set equipped with a partial operation $a*b$, will be said to be *(sub)additive* if $w(a*b) \leq w(a) + w(b)$ whenever $a*b$ is defined; and *strictly additive*, or *linear*, if $w(a*b) = w(a) + w(b)$ when this makes sense. The main property being the former, the prefix 'sub' will generally be omitted: for instance, an 'additively weighted' category has an 'additive' weight function on morphisms (3.1), in the first sense.

Occasionally, we shall also use the same category $\mathbf{w} = ([0, \infty], \geq)$ with the *multiplicative structure* $\mathbf{w}^\bullet$, where the tensor product is $\lambda.\mu$, and $\lambda.\infty = \infty$ for all $\lambda$ (cf. [G8]). Also here, a *multiplicative* function is actually submultiplicative.

**1.2. Directed metrics.** Now, a generalised metric space X, in the sense of Lawvere [Lw], called here a *directed metric space* or *$\delta$-metric space*, is a set X equipped with a $\delta$-*metric* $\delta$: X×X $\to [0, \infty]$, satisfying the axioms

(1) $\quad \delta(x, x) = 0,$ $\qquad\qquad\qquad\qquad \delta(x, y) + \delta(y, z) \geq \delta(x, z).$

This amounts to a category enriched over the symmetric monoidal closed category $\mathbf{w}^+$ considered above, with $\delta(x, y) = X(x, y)$ the hom-object in $[0, \infty]$. (If the value $\infty$ is forbidden, $\delta$ is often called a *quasi-pseudo-metric*, cf. [Ke]; but including it has crucial advantages, e.g. with respect to limits and colimits.)

$\delta$**Mtr** will denote the category of such $\delta$-metric spaces, with (weak) *contractions*, or *1-Lipschitz maps* f: X $\to$ Y, satisfying $\delta(x, x') \geq \delta(f(x), f(x'))$ for all x, x' $\in$ X, also called $\delta$-*maps*. Isomorphisms in this category are isometric - and will be called *isometric isomorphisms* or *1-Lipschitz isomorphisms* when needed. Limits and colimits exist and are calculated as in **Set**.

Products have the $l_\infty$-type $\delta$-metric (always defined because $\infty$ is included):

(2) $\quad \prod_i X_i, \quad \delta(\mathbf{x}, \mathbf{y}) = \sup \delta_i(x_i, y_i)$ $\qquad\qquad\qquad\qquad (\mathbf{x} = (x_i), \mathbf{y} = (y_i)).$

Equalisers have the restricted δ-metric. Sums have the obvious δ-metric (using ∞ also in the binary case):

(3)  $\Sigma_i X_i$,    $\delta((x, i), (y, i)) = \delta_i(x, y)$,    $\delta((x, i), (y, j)) = \infty$  $(i \neq j)$.

Coequalisers have the δ-metric induced on the quotient:

(4)  X/R,    $\delta(\xi, \eta) = \inf_{\mathbf{x}}(\Sigma_j \delta(x_{2j-1}, x_{2j}))$    $(\mathbf{x} = (x_1,..., x_{2p}); x_1 \in \xi; x_{2j} R x_{2j+1}; x_{2p} \in \eta)$.

The term 'directed' (used here) refers to the non-symmetric character of δ-metric spaces. And indeed such an object has a canonical preorder (to be used later, in 1.9)

(5)  $x \prec_\infty y$   if    $\delta(x, y) < \infty$.

More formally, a preordered set is the same as a δ-metric space with a truth-valued metric δ: X×X → **v**, taking values in {0, ∞}; the canonical preorder (5) gives thus the left adjoint of the embedding of preordered sets into δ-metric spaces (cf. 1.1.4).

The *reflected*, or *opposite*, δ-metric space $R(X) = X^{op}$ has the opposite δ-metric, $\delta^{op}(x, y) = \delta(y, x)$. A *symmetric* δ-metric, with $\delta = \delta^{op}$, is the same as an *écart* in Bourbaki [Bo]. More generally, a δ-metric space is *reflexive*, or self-dual, if it is isometrically isomorphic to its opposite (cf. 1.5). The notation $X \leq X'$ means that these δ-metric spaces have the same underlying set and $\delta_X \leq \delta_{X'}$, or equivalently that the identity of the underlying set is a δ-map $X' \to X$.

**1.3. Lipschitz maps.** We will also use the wider category $\delta_\infty\mathbf{Mtr}$ of δ-metric spaces and all *Lipschitz maps* f: X → Y, i.e. those mappings between the underlying sets for which the Lipschitz weight ‖f‖ is finite:

(1)  ‖f‖ = min $\{\lambda \in [0, \infty] \mid$ for all  $x, x' \in X, \delta(f(x), f(x')) \leq \lambda.\delta(x, x')\}$,

also called $\delta_\infty$-*maps*. A *Lipschitz isomorphism* will be an isomorphism of this category.

This category is finitely complete and cocomplete, with the same finite limits and colimits as $\delta\mathbf{Mtr}$; but note that, now, the δ-metric of a (co)limit is only determined up to Lipschitz-equivalence. Moreover, $\delta_\infty\mathbf{Mtr}$ is *multiplicatively weighted* [G8], with the Lipschitz weight: all identities have ‖$1_X$‖ ≤ 1 and composition gives ‖gf‖ ≤ ‖f‖.‖g‖. This also holds for $\delta\mathbf{Mtr}$, where ‖f‖ ≤ 1.

If X is a δ-metric space and $\lambda \in [0, \infty[$, we will write λX the same set equipped with the δ-metric $\lambda.\delta_X$ (where it is assumed that $\lambda.\infty = \infty$ for all λ, cf. [G8]). Thus, a $\delta_\infty$-map f: X → Y with ‖f‖ ≤ λ is the same as a δ-map λX → Y. More generally, as in [Lw], one can define λX where λ: $[0, \infty] \to [0, \infty]$ is any increasing mapping with $\lambda(0) = 0$ and $\lambda(\mu + \nu) \leq \lambda(\mu) + \lambda(\nu)$ (a lax monoidal functor $\mathbf{w}^+ \to \mathbf{w}^+$). For instance, the square-root mapping gives the δ-metric space √X.

**1.4. Tensor product.** The category $\delta\mathbf{Mtr}$ has a 'natural' symmetric monoidal closed structure (cf. [Lw], p. 153). The tensor product X⊗Y is the cartesian product of the underlying set, with the $l_1$-type δ-metric (instead of the $l_\infty$-type δ-metric of the categorical product)

(1)  $\delta((x, y), (x', y')) = \delta(x, x') + \delta(y, y')$;

it solves the usual universal problem, with respect to mappings which are 1-Lipschitz in each variable.

The exponential $Z^Y$ is the set of 1-Lipschitz maps Y → Z equipped with the δ-metric *of uniform convergence*

(2)  $\delta(h, k) = \sup_y \delta_Z(h(y), k(y))$ (with y varying in Y),

  $= \sup_{yy'} (\delta_Z(h(y), k(y')) - \delta_Y(y, y'))$ (with y, y' varying in Y).

The proof of the adjunction is standard (and can be deduced from the proof of theorem 5.6). The cartesian and tensor product are linked by the inequality

(3)  $X \times Y \leq X \otimes Y \leq 2.(X \times Y)$.

In $\delta_\infty\mathbf{Mtr}$, these products are isomorphic and denote isomorphic functors (in two variables). But we shall keep distinguishing such objects (and functors): the notation $X \times Y$ (resp. $X \otimes Y$) still denotes the *realisation* of the cartesian product given by the $l_\infty$-type (resp. $l_1$-type) $\delta$-metric.

**1.5. Standard models.** The line $\mathbf{R}$ and the standard interval $\mathbf{I}$ will have the euclidean metric $|x - y|$. Then, $\mathbf{R}^n$ and $\mathbf{I}^n$ have the *product metric*, $\sup_i |x_i - y_i|$, while $\mathbf{R}^{\otimes n}$ and $\mathbf{I}^{\otimes n}$ have the *tensor product metric*, $\Sigma_i |x_i - y_i|$.

But we are more interested in the following *non-symmetric* $\delta$-metrics. The *standard $\delta$-line* $\delta\mathbf{R}$ has the $\delta$-metric

(1)  $\delta(x, y) = y - x$, if $x \leq y$,  $\delta(x, y) = \infty$, otherwise;

its associated preorder is the natural order $x \leq y$ (1.2.5).

The *standard $\delta$-interval* $\delta\mathbf{I} = \delta[0, 1]$ has the subspace structure of the $\delta$-line. This also provides the cartesian powers $\delta\mathbf{R}^n$, $\delta\mathbf{I}^n$ and the tensor powers $\delta\mathbf{R}^{\otimes n}$, $\delta\mathbf{I}^{\otimes n}$. These $\delta$-metric spaces are not symmetric (for $n > 0$), but reflexive; in particular, the canonical reflecting isomorphism

(2)  $r: \delta\mathbf{I} \to (\delta\mathbf{I})^{op}$,  $t \mapsto 1 - t$,

will play a role, in *reflecting* paths and homotopies (in the opposite space).

The *standard $\delta$-circle* $\delta\mathbf{S}^1$ will be the coequaliser in $\delta\mathbf{Mtr}$ of the following two pairs of maps (equivalently)

(3)  $\partial^-, \partial^+: \{*\} \rightrightarrows \delta\mathbf{I}$,  $\partial^-(*) = 0$, $\partial^+(*) = 1$,

(4)  id, f: $\delta\mathbf{R} \rightrightarrows \delta\mathbf{R}$,  $f(x) = x + 1$.

Thus, the 'standard realisation' of the first coequaliser is the quotient $(\delta\mathbf{I})/\partial\mathbf{I}$, which identifies the endpoints; $\delta(x, y)$ takes values in $[0, 1[$, and can be viewed as measuring the length of the 'counterclockwise arc' from x to y, with respect to the whole circle. The structure $2\pi.\delta\mathbf{S}^1$ is also of interest (now, arcs are measured with respect to the radius).

More generally, the *n-dimensional $\delta$-sphere* will be the quotient of the monoidal $\delta$-cube $\delta\mathbf{I}^{\otimes n}$ modulo its (ordinary) boundary $\partial\mathbf{I}^n$,

(5)  $\delta\mathbf{S}^n = (\delta\mathbf{I}^n)/(\partial\mathbf{I}^n)$  $(n > 0)$,

while $\delta\mathbf{S}^0 = \{-1, 1\}$ will be given the *discrete* $\delta$-metric, infinite out of the diagonal (so that every mapping from this space to any other be a contraction). All $\delta$-spheres are reflexive.

**1.6. The symmetric case.** The full subcategory $\mathbf{Mtr}$ of symmetric $\delta$-metric spaces (1.2) has the same limits and colimits. Actually, it is reflective and coreflective in $\delta\mathbf{Mtr}$.

The coreflector, right adjoint to the embedding, is the well-known symmetrising procedure $D(x, x') = \delta(x, x') \vee \delta(x', x)$, based on the least symmetric δ-metric $D \geq \delta$. It will not be used here, since (for instance) it turns the δ-metric of δ**R** into the codiscrete δ-metric - infinite out of the diagonal.

But we shall frequently use the reflector $!(X, \delta) = (X, !\delta)$, based on the greatest symmetric δ-metric $!\delta \leq \delta$

(1) $\quad !\delta(x, x') = \inf_{\mathbf{x}}(\Sigma_j (\delta(x_{j-1}, x_j) \wedge \delta(x_j, x_{j-1}))) \qquad (\mathbf{x} = (x_0,..., x_p),\ x_0 = x,\ x_p = x')$,

which will be called the *symmetrised* δ-metric (of δ). The associated topology will be called the *symmetric topology* of the δ-metric space X, and is the one we are interested in (see 1.9).

This procedure turns the δ-metric of δ**R** into the euclidean metric. On δ**R**$^n$ the reflector gives a δ-metric $!(\delta\mathbf{R}^n)$ with ε-ball as in the second figure below, the convex hull of $[-\varepsilon, 0]^n \cup [0, \varepsilon]^n$

(2)
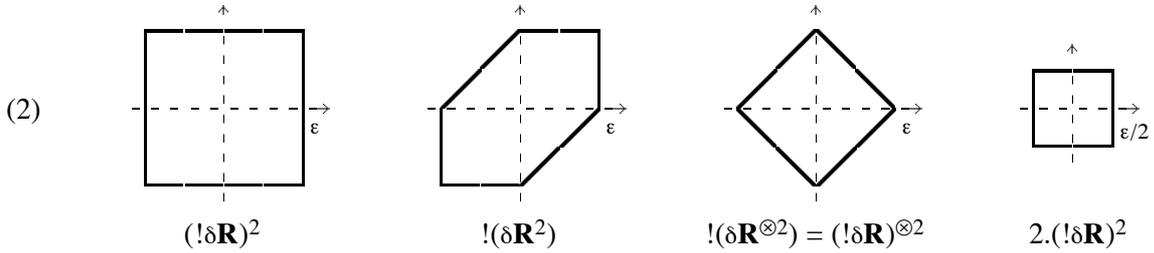

$\quad\quad (!\delta\mathbf{R})^2 \qquad\qquad !(\delta\mathbf{R}^2) \qquad\qquad !(\delta\mathbf{R}^{\otimes 2}) = (!\delta\mathbf{R})^{\otimes 2} \qquad 2.(!\delta\mathbf{R})^2$

while on the tensor powers $\delta\mathbf{R}^{\otimes n}$ it gives precisely the $l_\infty$-metric $(!\delta\mathbf{R})^{\otimes n}$, with ε-ball as above, in the third figure. All these δ-metrics are Lipschitz-equivalent, as follows from the figures above and from the following more general result.

**1.7. Proposition** (Symmetrisation and products). Given a finite family of δ-metric spaces $X_1,..., X_n$ we have the following inequalities (or equalities) for the δ-metrics obtained by symmetrisation, product and tensor product (on the cartesian product of the underlying sets $|X_1| \times ... \times |X_n|$)

(1) $\quad \Pi(!X_i) \leq !(\Pi X_i) \leq !(\otimes X_i) = \otimes(!X_i) \leq n.\Pi(!X_i)$,

so that all these δ-metrics are Lipschitz-equivalent and induce the same topology.

**Proof.** Recall the notation $\delta_X(x, x') = X(x, x')$, which comes from viewing a δ-metric space as an enriched category (1.2). The only non-standard point is the 'backward' inequality for the tensor product, proved in (c).

(a) First, to compare $\Pi(!X_i)$ and $!(\Pi X_i)$, note that

(2) $\quad \Pi(!X_i)(\mathbf{x}, \mathbf{y}) = \sup (!X_i)(x_i, y_i) \leq \sup X_i(x_i, y_i) = (\Pi X_i)(\mathbf{x}, \mathbf{y})$;

since $\Pi(!X_i)$ is symmetric, it follows that $\Pi(!X_i) \leq !(\Pi X_i)$.

(b) The second inequality, $!(\Pi X_i) \leq !(\otimes X_i)$, is a straightforward consequence of $(\Pi X_i) \leq (\otimes X_i)$.

(c) We prove now that $!(\otimes X_i) = \otimes(!X_i)$. The δ-metric of the latter is:

(3) $\quad \otimes(!X_i)(\mathbf{x}, \mathbf{y}) = \Sigma_i (!X_i)(x_i, y_i) \leq \Sigma_i X_i(x_i, y_i) = (\otimes X_i)(\mathbf{x}, \mathbf{y})$.

Since $\otimes(!X_i)$ is symmetric, we have $\otimes(!X_i) \leq !(\otimes X_i)$. The opposite inequality is more subtle: take a sequence of n+1 points $\mathbf{z}^j$, which varies from $\mathbf{x}$ to $\mathbf{y}$, by changing one coordinate at a time

(4)  $\mathbf{z}^j = (y_1,..., y_j, x_{j+1},..., x_n)$, $\qquad \mathbf{z}^0 = \mathbf{x}$, $\mathbf{z}^n = \mathbf{y}$  $\qquad (j = 0,..., n)$.

and apply the triangular inequality:

(5)  $!(\otimes X_i)(\mathbf{x}, \mathbf{y}) \leq \Sigma_j !(\otimes X_i)(\mathbf{z}^{j-1}, \mathbf{z}^j)$.

Now, $\mathbf{z}^{j-1}$ and $\mathbf{z}^j$ only differ at the j-th coordinate ($x_j$ or $y_j$, respectively); restricting the domain of the 'inf' in the right term above to those sequences in $!(\otimes X_i)$ where only the j-th coordinate changes, we get the δ-metric $!X_j$, and the inequality:

(6)  $\Sigma_j !(\otimes X_i)(\mathbf{z}^{j-1}, \mathbf{z}^j) \leq \Sigma_j (!X_j)(x_j, y_j) = \otimes(!X_i)(\mathbf{x}, \mathbf{y})$.

(d) Finally, the last inequality in (1) is obvious. □

**1.8. Definition and Proposition** (The length of paths). Let X be a δ-metric space and a: $\mathbf{I} \to X$ a set-theoretical mapping. We define its *span* sp(a) and its *length* L(a) with the following functions, taking values in $[0, \infty]$:

(1)  $sp(a) = \sup_{\mathbf{t}} \delta(a(t_0), a(t_1))$  $\qquad (0 \leq t_0 < t_1 \leq 1)$,

$L(a) = \sup_{\mathbf{t}} L_{\mathbf{t}}(a)$,

$L_{\mathbf{t}}(a) = \Sigma_i \delta(a(t_{i-1}), a(t_i))$  $\qquad (0 = t_0 < t_1 < ... < t_p = 1)$.

These functions satisfy the following properties, where $0_x$ is the constant paths at (any) point x, a + b denotes (any) path-concatenation of consecutive paths and $\|a\|$ is the Lipschitz weight (1.3)

(a)  $sp(0_x) = L(0_x) = 0$,

(b)  $sp(a+b) \leq sp(a) + sp(b)$,  $\qquad L(a+b) = L(a) + L(b)$,

(c)  $sp(a\rho) \leq sp(a)$,  $\qquad L(a\rho) \leq L(a)$  $\qquad$ (for every *weakly increasing map* $\rho: \mathbf{I} \to \mathbf{I}$),

(d)  $sp(a\rho) = sp(a)$,  $\qquad L(a\rho) = L(a)$  $\qquad$ (for every *increasing homeomorphism* $\rho: \mathbf{I} \to \mathbf{I}$),

(e)  $sp(a, b) = sp(a) + sp(b)$,  $\qquad L(a, b) = L(a) + L(b)$  $\qquad$ (for all paths $(a, b): \mathbf{I} \to X \otimes Y$),

(f)  $sp(a) \leq L(a) \leq \|a\|$,

(g)  L is the least function on set-theoretical paths which is strictly additive for concatenation, invariant for reparametrisation on increasing homeomorphisms $\mathbf{I} \to \mathbf{I}$, and satisfies $L \geq sp$,

(h)  $sp(f \circ a) \leq \|f\|. sp(a)$,  $\qquad L(f \circ a) \leq \|f\|. L(a)$  $\qquad$ (for all $\delta_\infty$-maps $f: X \to Y$),

(i)  for a: $\mathbf{I} \to \delta\mathbf{R}^{\otimes n}$, we get $L(a) = sp(a) = \Sigma_i (a_i(1) - a_i(0))$ if a is (weakly) increasing, and $L(a) = \infty$ otherwise.

Finally, note that the length L(a) can be finite even when a is not Lipschitz: $\|a\| = \infty$.

**Proof.** The properties of the span being obvious, we only verify the ones of the length. It will be useful to note that, if the partition $\mathbf{t}'$ is finer than $\mathbf{t}$, then $L_{\mathbf{t}}(a) \leq L_{\mathbf{t}'}(a)$, because of the triangular property of δ-metrics.

Point (a) is obvious. For (b), the inequality $L(a+b) \leq L(a) + L(b)$ follows easily: given a partition $\mathbf{t}$ for $c = a + b$, call $\mathbf{t}'$ its refinement by introducing the point $t = 1/2$ (if missing); thus $L_{\mathbf{t}}(c) \leq L_{\mathbf{t}'}(c)$ and the latter term can be split into two summands $\leq L(a) + L(b)$. For the other inequality, it is sufficient to note that a partition for a and one for b determine a partition of [0, 2], which can be scaled down to the standard interval.

Points (c) is obvious, since $L(a\rho)$ is computed on the partitions of $\rho[0, 1] \subset [0, 1]$; (d) is a consequence. For (e), the inequality $L(a, b) \leq L(a) + L(b)$ is obvious, and the other follows again from the first remark: given a partition for $a$ and one for $b$, by using a common refinement for both we get higher values.

Finally, (f) - (h) are plain; (i) is obvious for $n = 1$, and follows from (e) for higher $n$. For the last remark, taking $X = \delta\mathbf{R}$, the square-root map $f: \mathbf{I} \to \mathbf{R}$ is not Lipschitz but has a finite length, as any increasing path: $L(f) = f(1) - f(0) = 1$. □

**1.9. The associated topology and direction.** A $\delta$-metric space $X$ has an obvious *directed topology*, generated by the *open discs* $\{x \mid \delta(x_0, x) < \varepsilon\}$. *We will not use this construct.*

Rather, we shall use the *symmetric topology*, determined as above by the symmetric $\delta$-metric $!\delta$ defined above (1.6.1). And we will keep a trace of the 'directed' information of the original $\delta$, in various ways, more or less effective, based on the settings for directed homotopy used in [G3].

The simplest, if rather poor way, is by the *associated preorder* $x \prec_\infty y$, defined by $\delta(x, y) < \infty$ (1.2.5). We have thus a forgetful functor with values in the category p**Top** of preordered spaces and continuous preorder-preserving mappings

(1)   **p**: $\delta$**Mtr** $\to$ p**Top**,              ($\delta_\infty$**Mtr** $\to$ p**Top**),

which preserves finite products (since the symmetrising functor $!: \delta$**Mtr** $\to$ **Mtr** preserves finite products up to Lipschitz equivalence (1.7), and $\delta(x, y) < \infty$ if and only if this holds on all components.

Thus, $\mathbf{p}(\delta\mathbf{R}^n) = \mathbf{p}(\delta\mathbf{R}^{\otimes n}) = \uparrow\mathbf{R}^n$, the euclidean n-dimensional space with the product order. Similarly, $\mathbf{p}(\delta\mathbf{I}^n) = \mathbf{p}(\delta\mathbf{I}^{\otimes n}) = \uparrow\mathbf{I}^n$. But a preorder is a poor way of describing direction, which does not allow for non-reversible loops. Thus, $\mathbf{p}(\delta\mathbf{S}^1)$ gets the chaotic preorder and misses any information of direction.

A more accurate way of keeping the 'directed' information of the original $\delta$ is using *d-spaces*, or *spaces with distinguished paths* [G3]. We have thus a forgetful functor

(2)   **d**: $\delta$**Mtr** $\to$ d**Top**,              ($\delta_\infty$**Mtr** $\to$ d**Top**),

which equips a $\delta$-metric space $X$ with the associated symmetric topology *and* the d-structure where a (continuous) path $a: \mathbf{I} \to X$ is distinguished if and only if it is L-*feasible*, i.e. it has a finite length $L(a)$ (1.8.1). The axioms of d-spaces are satisfied (by 1.8): distinguished paths contain all the constant ones, are closed under concatenation and partial reparametrisation by weakly increasing maps $\mathbf{I} \to \mathbf{I}$. And, of course, a Lipschitz map $f: X \to Y$ of $\delta$-metric spaces preserves feasible paths. It will be relevant to note that this functor takes the tensor (or cartesian) product in $\delta_\infty$**Mtr** (or $\delta$**Mtr**) to the cartesian product of d-spaces (where a path is distinguished if and only if its two components are).

Now, in d**Top**, $\mathbf{d}(\delta\mathbf{S}^1) = \uparrow\mathbf{S}^1 = \uparrow\mathbf{I}/\partial\mathbf{I}$ is the standard directed circle [G3], where path are only allowed to turn in a given direction. But the functor **d** need not preserve quotients: for instance, $\mathbf{d}(\delta\mathbf{R}) = \uparrow\mathbf{R}$ is the standard directed line, with the increasing paths as distinguished ones; the quotient d-space $\uparrow\mathbf{R}/G_\vartheta$, modulo the action of the dense subgroup $G_\vartheta = \mathbf{Z} + \vartheta\mathbf{Z}$ (for an irrational number $\vartheta$), is a non-trivial object, with the homology of the 2-dimensional torus (cf. [G6]), while $(\delta\mathbf{R})/G_\vartheta$ has a trivial $\delta$-metric, always zero. A finer notion of 'weighted space', studied in Sections 5-7, will be able to express such phenomena within *weighted* algebraic topology (not just the *directed* one).

Note that the forgetful functor d**Top** → p**Top** provided by the path-preorder $x \preceq x'$ (there is a distinguished path from $x$ to $x'$), applied to a d-space of type **d**X, gives a *finer* preorder than **p**X, generally more interesting than the latter (two points at a finite distance may be disconnected, or not linked by a feasible path.)

One could also try an intermediate way of codifying direction, by some notion of *local preorder*, as in Krishnan [Kr]. For instance, one could use a family of preorders, indexed on the open subsets $U$ of $X$ (always with respect to the symmetric topology)

(3)   $x \prec_U x'$  if  ($x \preceq x'$ in the path-preorder of $U$, with respect to the induced δ-metric),

so that, if $U \subset V$, $x \prec_U x'$ implies $x \prec_V x'$. Thus, in $\delta S^1$, every proper open arc $U$ gets a total order (Lipschitz isomorphic to an open interval).

## 2. Elementary and extended homotopies

The standard δ-interval δ**I** generates a cylinder endofunctor, which yields *elementary* homotopies in δ**Mtr**, and *extended* homotopies in $\delta_\infty$**Mtr**. We need both, but only the latter can be concatenated. The letter α denotes an element of the set {0, 1}, written −, + in superscripts.

**2.1. Elementary and extended paths.** Let $X$ be a δ-metric space. An *elementary path* (resp. an *extended path*, or *Lipschitz path*) in $X$ will be a 1-Lipschitz (resp. a Lipschitz) map a: δ**I** → X.

Thus, a set-theoretical mapping a: **I** → X is an elementary path if and only if $\|a\| \leq 1$, for the Lipschitz weight (1.3.1)

(1)   $\|a\| = \min \{\lambda \in [0, \infty] \mid \text{for all } t \leq t' \text{ in } [0, 1], \delta(a(t), a(t')) \leq \lambda.(t' - t)\}$,

and is an extended path if and only if $\|a\| < \infty$. Elementary paths cannot be concatenated, because - loosely speaking - this procedure doubles the velocity, whose least upper bound is the Lipschitz weight. Recall that the (finite) length $L(a) \leq \|a\|$ has been defined in 1.8.

The *reflected* (elementary or extended) path is obtained in the obvious way

(2)   $a^{op} = ar: \delta\mathbf{I} \to X^{op}$,                    $r(t) = 1 - t$,

A *reversible* extended path is a mapping a: **I** → X such that both a and $a^{op}$ are extended paths δ**I** → X. This amounts to a Lipschitz map a: !δ**I** → X, with respect to the ordinary metric $|t - t'|$ of the euclidean interval.

**2.2. Concatenation.** Given two consecutive Lipschitz paths a, b: δ**I** → X, with a(1) = b(0), the usual construction gives a concatenated path

(1)   a + b: δ**I** → X,                    $\| a + b \| \leq 2.(\|a\| \vee \|b\|)$,

(which need not be elementary when a and b are). Let us review this fact in a more formal way.

In the category δ**Mtr**, pasting two copies of the standard δ-interval, one after the other, can be realised as δ[0, 2] ⊂ δ**R**, or (isometrically) as 2.δ**I** (with the double δ-metric)

(2)
$$\begin{array}{ccc} \{*\} & \xrightarrow{\partial^+} & \delta\mathbf{I} \\ \partial^- \downarrow & \phantom{--} \downarrow k^- \\ \delta\mathbf{I} & \xrightarrow{k^+} & 2.\delta\mathbf{I} \end{array} \qquad k^-(t) = t/2, \quad k^+(t) = (t+1)/2.$$

Of course, this is of no help to concatenate *elementary* paths. *Moving to the category* $\delta_\infty\mathbf{Mtr}$, this pushout can be realised as the Lipschitz-isomorphic object $\delta\mathbf{I}$. This yields the *standard concatenation pushout* (the left diagram below, which actually lives in $\delta\mathbf{Mtr}$ but is not a pushout there)

(3)
$$\begin{array}{ccc} \{*\} & \xrightarrow{\partial^+} & \delta\mathbf{I} \\ \partial^- \downarrow & \phantom{--} \downarrow k^- \\ \delta\mathbf{I} & \xrightarrow{k^+} & \delta\mathbf{I} \end{array} \qquad \begin{array}{ccc} X & \xrightarrow{\partial^+} & X\otimes\delta\mathbf{I} \\ \partial^- \downarrow & \phantom{--} \downarrow k^- \\ X\otimes\delta\mathbf{I} & \xrightarrow{k^+} & X\otimes\delta\mathbf{I} \end{array} \qquad \begin{array}{l} k^-(x, t) = (x, t/2), \\ k^+(x, t) = (x, t+1/2). \end{array}$$

*This pushout is preserved by any functor* $X\otimes-$ (yielding the right hand pushout above), or by $X\times-$. In fact, $X\otimes-: \delta\mathbf{Mtr} \to \delta\mathbf{Mtr}$ preserves the pushout (2), as a left adjoint; and the embedding $\delta\mathbf{Mtr} \to \delta_\infty\mathbf{Mtr}$ preserves finite colimits (1.3).

**2.3. The elementary cylinder.** The $\delta$-interval $\delta\mathbf{I}$ is an internal lattice *in* $(\delta\mathbf{Mtr}, \otimes)$: its structure consists of two *faces* $(\partial^-, \partial^+)$, a *degeneracy* (e), two *connections* or main operations $(g^-, g^+)$ and an *interchange* (s)

(1) $\quad \{*\} \underset{e}{\overset{\partial^\alpha}{\rightrightarrows}} \delta\mathbf{I} \overset{g^\alpha}{\Leftarrow} \delta\mathbf{I}^{\otimes 2} \qquad\qquad s: \delta\mathbf{I}^{\otimes 2} \to \delta\mathbf{I}^{\otimes 2} \qquad\qquad (\alpha = \pm)$

$$\partial^\alpha(*) = \alpha, \qquad g^-(t, t') = t \vee t', \qquad g^+(t, t') = t \wedge t', \qquad s(t, t') = (t', t).$$

As a consequence, the *elementary cylinder* endofunctor

(2) $\quad \mathrm{I}: \delta\mathbf{Mtr} \to \delta\mathbf{Mtr}, \qquad\qquad \mathrm{I}(-) = - \otimes \delta\mathbf{I},$

has natural transformations, which will be denoted by the same symbols and names

(3) $\quad 1 \underset{e}{\overset{\partial^\alpha}{\rightrightarrows}} \mathrm{I} \overset{g^\alpha}{\Leftarrow} \mathrm{I}^2 \qquad\qquad s: \mathrm{I}^2 \to \mathrm{I}^2.$

This functor $\mathrm{I}$ and these transformations satisfy the axioms of a *cubical monad with interchange* [G1, G2]:

(4) $\quad e\partial^\alpha = 1, \qquad\qquad eg^\alpha = e.\mathrm{I}e \;(= e.e\mathrm{I}) \qquad\qquad$ (*degeneracy axiom*),

$\phantom{(4)}\quad g^\alpha.\mathrm{I}g^\alpha = g^\alpha.g^\alpha\mathrm{I}, \qquad g^\alpha.\mathrm{I}\partial^\alpha = 1 = g^\alpha.\partial^\alpha\mathrm{I} \qquad\qquad$ (*associativity, unit*),

$\phantom{(4)}\quad g^\beta.\mathrm{I}\partial^\alpha = \partial^\alpha e = g^\beta.\partial^\alpha\mathrm{I} \qquad\qquad\qquad\qquad\qquad$ (*absorbency*, for $\alpha \ne \beta$),

$\phantom{(4)}\quad s.s = 1, \qquad\qquad \mathrm{I}e.s = e\mathrm{I},$

$\phantom{(4)}\quad s.\mathrm{I}\partial^\alpha = \partial^\alpha\mathrm{I}, \qquad\qquad g^\alpha.s = g^\alpha \qquad\qquad\qquad\qquad$ (*interchange*).

The cylinder I has a *generalised reversion*, via the *reflection* of δ-metric spaces (as always happens in *directed* algebraic topology, e.g. for d-spaces [G3] and for differential graded algebras [G1])

(5)  $rX = X \otimes r: I.RX \to R.IX,$   $(x, t) \mapsto (x, 1 - t),$

  $RrR.r = id,$   $Re.r = eR,$   $r.\partial^- R = R\partial^+,$

  $r.g^- R = Rg^+.r_2,$   $Rs.r_2 = r_2.sR.$

where $r_2 = rI.Ir: (I^2R \to IRI \to RI^2)$ is the reversion of the double cylinder.

Within δ**Mtr**, the cylinder endofunctor $I = - \otimes \delta I$ has a right adjoint, the *elementary-path functor*, or *elementary cocylinder*

(6)  P: δ**Mtr** → δ**Mtr**,   $P(Y) = Y^{\delta I}.$

The δ-metric space $Y^{\delta I}$ is the set of elementary paths δ**Mtr**(δI, Y) with the δ-metric of uniform convergence (1.4.2). The lattice structure of δI in d**Top** produces – contravariantly – a dual structure on P (a cubical *co*monad with interchange [G1, G2]).

**2.4. Homotopies.** An *elementary homotopy* φ: f → g: X → Y is defined as a δ-map φ: IX = X⊗δI → Y whose two faces are f and g, respectively: $\partial^-(\varphi) = \varphi \circ \partial^- = f$, $\partial^+(\varphi) = \varphi \circ \partial^+ = g$. In particular, an elementary path is a homotopy between two points, a: x → x': {∗} → X.

More generally, an *extended homotopy*, or *Lipschitz homotopy*, is a Lipschitz map φ: X×δI → Y; and an extended path is an extended homotopy between two points, a: x → x': {∗} → X. (Note that a Lipschitz map defined on the singleton is always a δ-map.) An extended homotopy has a Lipschitz weight ‖φ‖, which is ≤ 1 if and only if φ is elementary.

In both cases, the main operations produced by the cylinder functor (for φ: f → g: X → Y; u: X' → X; v: Y → Y'; ψ: g → h: X → Y) are:

(a) *whisker composition* of (elementary or extended) maps and homotopies

  v∘φ∘u: vfu → vgu     (v∘φ∘u = v.φ.Iu: IX' → Y'),

(b) *trivial homotopies*:   $0_f$: f → f     ($0_f$ = fe: IX → Y),

and satisfy obvious axioms for units and associativity:

(1)  $1 \circ \varphi \circ 1 = \varphi,$  $v \circ 0_f \circ u = 0_{vfu},$  $(v'v) \circ \varphi \circ (uu') = v' \circ (v \circ \varphi \circ u) \circ u'.$

This 2-dimensional structure, weaker than a bicategory or a sesquicategory (there is no vertical composition of 2-cells) is called an *h-category* in [G2] (it is a category enriched over reflexive graphs, with a suitable monoidal structure).

More precisely, we will write:

(i)  $\delta_1$**Mtr**  the h-category of δ-metric spaces, with *weak contractions* and *elementary homotopies*,

(ii)  $\delta_\infty$**Mtr**  the h-category of δ-metric spaces, with *Lipschitz maps* and *Lipschitz homotopies*,

(iii)  δ**Mtr**  the intermediate h-category of *weak contractions* and *Lipschitz homotopies* between them.

Actually, $\delta_1$**Mtr** is a bicomplete IP-*homotopical category* [G1]: it has adjoint functors I ⊣ P, with the required structure (a cubical monad and a cubical comonad, respectively), colimits (preserved by the cylinder) and limits (preserved by the cocylinder). Therefore, all results of [G1] for such a structure apply (as for cochain algebras).

Moreover, in $\delta_\infty\mathbf{Mtr}$ (and $\delta\mathbf{Mtr}$) consecutive homotopies can be pasted via the *concatenation pushout* of the cylinder functor (the right-hand diagram in 2.2.3). The concatenation $\varphi+\psi$ of two consecutive homotopies $(\partial^+\varphi = \partial^-\psi)$ is thus computed as usual:

(2) $\quad (\varphi+\psi)(x, t) = \varphi(x, 2t)\quad$ for $0 \leq t \leq 1/2,\quad = \psi(x, 2t-1)\quad$ for $1/2 \leq t \leq 1$.

The *reflected* homotopies (elementary or extended) live in the opposite 'spaces' (as for paths, 2.1.2)

(3) $\quad \varphi^{op}: Rg \to Rf: RX \to RY, \qquad\qquad \varphi^{op} = R\varphi.rX = (\delta\mathbf{I}(RX) \to R(\delta\mathbf{I}X) \to RY)$.

As always in directed algebraic topology, *homotopy equivalence is a complex notion*, which has to be considered not only for 'spaces' but also for their algebraic counterpart - weighted categories. This will be briefly considered in Section 5. A more complete study for d-spaces can be found in [G3].

**2.5. Double homotopies and 2-homotopies.** Extended homotopies in dimension 2, based on the *second order cylinder* $\mathbf{I}^2 X = X \otimes \delta\mathbf{I}^2$, can now be treated like for d-spaces, in [G3]. Roughly speaking, they behave as in **Top**, as long as we work on the standard square $[0, 1]^2$ with *increasing* Lipschitz maps. We will only sketch the main points; the interested reader can look at [G3], 2.5 - 2.6.

An extended *double homotopy* is a map $\Phi: X \otimes \delta\mathbf{I}^2 = \mathbf{I}^2 X \to Y$; it has four faces $\partial_i^\alpha(\Phi)$, with $\alpha = \pm$ and $i = 1, 2$. The concatenation of extended double homotopies *in direction* 1 or 2 is defined as usual (under the obvious boundary conditions) and satisfies a strict *middle-four interchange property*.

An extended *2-homotopy* $\Phi: \varphi \to \psi: f \to g: X \to Y$ is a double homotopy whose faces $\partial_1^\alpha$ are degenerate, while the faces $\partial_2^\alpha$ are $\varphi, \psi$ (the other choice is equivalent, by interchange). Such particular double homotopies are closed under concatenations in each direction (also because $0_f + 0_f = 0_f$). The preorder $\varphi \prec_2 \psi$ (i.e., there is a 2-homotopy $\varphi \to \psi$) spans an equivalence relation $\simeq_2$.

Various constructions of 2-dimensional homotopies are crucial tools in the theory.

(a) Two 'horizontally' consecutive extended homotopies

(1) $\quad \varphi: f^- \to f^+: X \to Y, \qquad\qquad \psi: g^- \to g^+: Y \to Z$,

can be composed, to form an extended double homotopy $\psi \circ \varphi$ (which is elementary if $\varphi$ and $\psi$ are)

(2)
$$\begin{array}{ccc} g^-f^- & \xrightarrow{g^-\circ\varphi} & g^-f^+ \\ \psi\circ f^- \downarrow & \psi\circ\varphi & \downarrow \psi\circ f^+ \\ g^+f^- & \xrightarrow{g^+\circ\varphi} & g^+f^+ \end{array} \qquad \begin{array}{l} \psi\circ\varphi = \psi.(\varphi\otimes\delta\mathbf{I}): X\otimes\delta\mathbf{I}^2 \to Y\otimes\delta\mathbf{I} \to Z, \\ \partial_1^\alpha(\Phi) = \psi.(\varphi\partial^\alpha\otimes\delta\mathbf{I}) = \psi\circ f^\alpha, \\ \partial_2^\alpha(\Phi) = \psi.(\varphi\otimes\partial^\alpha) = g^\alpha\circ\varphi. \end{array}$$

(Together with the whisker composition, in 2.4, this is a particular instance of the cubical enrichment produced by the (co)cylinder functor: composing a p-uple homotopy $\Phi: \mathbf{I}^p X \to Y$ with a q-uple one $\Psi: \mathbf{I}^q Y \to Z$ gives a (p+q)-uple homotopy $\Psi \circ \Phi = \Psi.\mathbf{I}^q\Phi$.)

(b) *Acceleration.* For every extended homotopy $\varphi: f \to g$, there are *acceleration* extended 2-homotopies (cf. [G3], 2.6.4)

(3) $\quad \Theta': 0_f + \varphi \to \varphi, \qquad\qquad \Theta'': \varphi \to \varphi + 0_g$,

(but *not* the other way round: slowing down conflicts with direction).

(c) *Folding*. A double extended homotopy $\Phi: A \otimes \delta \mathbf{I}^2 \to X$ with faces $\varphi, \psi, \sigma, \tau$ (as below) produces a 2-homotopy $\Psi$, by pasting $\Phi$ with two double homotopies produced by the connections $g^\pm$ (2.3), denoted with #

(4)
$$\begin{array}{ccccccc} f & = & f & \xrightarrow{\sigma} & h & \xrightarrow{\psi} & g \\ \| & \# \;\varphi \downarrow & & \Phi & \downarrow \psi \; \# & & \| \\ f & \xrightarrow{\varphi} & k & \xrightarrow{\tau} & g & = & g \end{array} \qquad \Psi: 0_f + \sigma + \psi \to \varphi + \tau + 0_g : f \to g.$$

Together with accelerations, this 2-homotopy $\Psi$ shows that the faces of $\Phi$ satisfy: $\sigma + \psi \simeq_2 \varphi + \tau$.

**2.6. Homotopy constructs.** We have now the tools to develop for $\delta$-metric spaces the homotopy constructs of [G4]: (a) homotopy pushouts and pullbacks; (b) mapping cones, suspension and cofibration sequences; (c) homotopy fibres, loop-objects and fibration sequences (in the pointed case).

Then, the higher properties of this machinery need concatenation, and can only be developed in $\delta_\infty \mathbf{Mtr}$. We will not write down these computations here.

# 3. The fundamental weighted category

The fundamental category of a $\delta$-metric space is defined. Non obvious computations are based on a van Kampen-type theorem (3.6), similar to R. Brown's version for the fundamental groupoid of spaces [Br] or to our version for spaces with distinguished paths [G3].

**3.1. Weighted categories.** An *additively weighted category* is a category $X$ where every morphism $a$ is equipped with a *weight*, or *cost*, $w(a) \in [0, \infty]$, so that two obvious axioms are satisfied, for identities and composition (written in additive notation):

($w^+$cat.0)  $w(0_x) = 0$, for all objects $x$ of $X$,

($w^+$cat.1)  $w(a + b) \leq w(a) + w(b)$, for all pairs of consecutive arrows $a, b$.

This is called a normed category in [Lw, BG]; see the Introduction and [G8] for the present terminology about the additive and multiplicative variants. We will omit the term 'additive' when there is no ambiguity; we also speak of a $w^+$-*category*, for short. The weight is said to be *linear*, or *strictly additive*, if $w(a + b) = w(a) + w(b)$ for all composites.

A $w^+$-*functor* $f: X \to Y$, or *1-Lipschitz functor*, is a functor between such categories which satisfies the condition $w(f(a)) \leq w(a)$, for all morphisms $a$ of X. A $w^+$-*transformation* $\varphi: f \to g$ is a natural transformation between w-functors. All this forms the 2-category $w^+\mathbf{Cat}$ of (small) $w^+$-categories, also written $w\mathbf{Cat}$. (Here we do not use the multiplicative analogue, for which see [G8]).

Weighted categories can be viewed as categories enriched over the symmetric monoidal closed category $w^+\mathbf{Set}$ of *weighted sets* (i.e., sets equipped with a mapping $w: X \to \mathbf{w^+}$, and maps $f: X \to Y$ such that $w(f(x)) \leq w(x)$, for all $x \in X$), with an 'additive' tensor product [BG, G8].

The *opposite weighted category* $X^{op}$ is the opposite category with the 'same' weight.

Also here (as for $\delta$-metric spaces, in 1.3) we need the wider category $w_\infty \mathbf{Cat}$ of weighted categories and *Lipschitz functors* f: X → Y, or $w_\infty$-*functors*, i.e. the functors between weighted categories having a finite Lipschitz weight

(1)  $\|f\| = \min \{\lambda \in [0, \infty] \mid \text{for all } a \in \text{Mor}(X),\ w_Y(f(a)) \leq \lambda.w_X(a)\}$.

With this weight, the category $w_\infty \mathbf{Cat}$ is *multiplicatively weighted* (cf. 1.3). This also holds for w**Cat**, which is the wide subcategory of the functors f such that $\|f\| \leq 1$.

**3.2. Proposition** (The monoidal closed structure). The category w**Cat** has a symmetric monoidal closed structure, with tensor product $X \otimes Y$ consisting of the cartesian product of the underlying categories with the $l_1$-type weight on a map (a, b): (x, y) → (x', y'):

(1)  $w_\otimes(a, b) = w(a) + w(b)$.

The internal hom $Z^Y$ is the category of 1-Lipschitz functors h: Y → Z and *all* natural transformations φ: h → k between such functors, with the (plainly subadditive) weight:

(2)  $W(\varphi) = \sup_y w_Z(\varphi(y))$,                                                      (y ∈ ObY).

**Proof.** First a 1-Lipschitz functor f: X⊗Y → Z defines a functor g: X → $Z^Y$, sending an object x to the 1-Lipschitz functor

(3)  $g(x) = f(x, -): Y \to Z$,           $w_Z(g(x)(b)) = w_Z(f(0_x, b)) \leq w_\otimes(0_x, b) = w(b)$,

and the X-morphism a: x → x' to the natural transformation g(a) = f(a, – ): g(x) → g(x'). The functor g itself is a contraction:

(4)  $W(g(a)) = \sup_y w_Z(g(a)(y)) = \sup_y w_Z(f(a, 0_y)) \leq \sup_y w_\otimes(a, 0_y) = w(a)$.

Conversely, given a 1-Lipschitz functor g: X → $Z^Y$, we define the functor f: X⊗Y → Z in the usual, obvious way, and verify that it is 1-Lipschitz, on a map (a, b): (x, y) → (x', y')

(5)  $w_Z(f(a, b)) = w_Z(g(a)(y) + g(x)(b)) \leq w_Z(g(a)(y)) + w_Z(g(x)(b)) \leq w(a) + w(b)$,

where the last inequality follows from:

(6)  $w_Z(g(a)(y)) \leq \sup_{y'} w_Z(g(a)(y')) = W(g(a)) \leq w(a)$,           $w_Z(g(x)(b)) \leq w(b)$.           □

**3.3. Homotopy for weighted categories.** This elementary theory is based on the *directed interval* **2** = {0 → 1}, an order category on two objects, where the only non-trivial arrow has weight w(0 → 1) = 1. The obvious *faces* $\partial^\pm$: **1** → **2** are defined on the singleton category **1** = {∗}.

A *point* x: **1** → X of a small weighted category X is an object of the latter; we will also write x ∈ X. An *extended path* a: **2** → X from x to x' amounts to a *feasible* arrow a: x → x' of X *with a finite weight* w(a); in fact, the latter coincides with the Lipschitz weight $\|a\|$, as a functor on **2**; concatenation of extended paths amounts to composition in X (strictly associative, with strict identities). Elementary paths amount to arrows with w(a) ≤ 1, which will also be called 1-*morphisms*, and cannot be concatenated, generally.

The monoidal closed structure described above (3.2) gives the *elementary cylinder* $IX = X \otimes \mathbf{2}$ and its right adjoint, the *elementary path* $PY = Y^{\mathbf{2}}$, which consists of the category of 1-morphisms of Y and their arbitrary commutative squares, with weight $W(b, b') = w(b) \vee w(b')$. This shows that an *elementary homotopy*, or *elementary natural transformation*, $\varphi: f \to g: X \to Y$ is the same as a natural transformation between 1-Lipschitz functors, *satisfying* $w(\varphi(x)) \leq 1$ for all $x \in X$. Such homotopies cannot be concatenated, since their vertical composition is no longer elementary, in general. An *elementary isomorphism* of 1-Lipschitz functors will be an elementary natural transformation having an inverse in the same domain; plainly, this amounts to an invertible natural transformation such that $w(\varphi(x)) = 1 = w(\varphi^{-1}(x))$ for all points x.

Let us consider now an arbitrary functor $\varphi: |X \otimes \mathbf{2}| \to |Y|$ *between the underlying categories*, and define its *reduced weight* $|\varphi|$ as follows (writing $\varphi(x) = \varphi(x, 0 \to 1)$, for $x \in X$)

(1) $\quad |\varphi| = \sup_x w_Y(\varphi(x)) = \min \{\lambda \in [0, \infty] \mid \text{for all } x \in X, \ w_Y(\varphi(x)) \leq \lambda\}$.

Note now that $\varphi$ sends a map $(a, 1_0)$ of $X \otimes \mathbf{2}$ to $f(a)$, a map $(a, 1_1)$ to $g(a)$ and a map $(a: x \to x', 0 \to 1)$ to

(2) $\quad b = \varphi(a, 0 \to 1) = g(a) \circ \varphi(x) = \varphi(x') \circ f(a)$,

$\quad w(b) \leq w(f(a)) + w(\varphi(x')) \leq \|f\|.w(a) + |\varphi| \leq (\|f\| \vee |\varphi|).(w(a) + 1)$,

$\quad w(b) \leq w(g(a)) + w(\varphi(x)) \leq \|g\|.w(a) + |\varphi| \leq (\|g\| \vee |\varphi|).(w(a) + 1)$.

Therefore, the Lipschitz weight of $\varphi$ as a functor $\varphi: |X \otimes \mathbf{2}| \to |Y|$, or *global weight*, is

(3) $\quad \|\varphi\| = \|f\| \vee \|g\| \vee |\varphi|$.

A *Lipschitz homotopy*, or *Lipschitz natural transformation*, $\varphi: f \to g: X \to Y$ will be a natural transformation with a finite global weight $\|\varphi\|$, or equivalently a natural transformation of *Lipschitz functors* with a finite reduced weight $|\varphi|$. We will write:

(i)   $w_1\mathbf{Cat}$ the h-category of weighted categories, with 1-Lipschitz functors and elementary natural transformations,

(ii)  $w_\infty\mathbf{Cat}$ the 2-category of weighted categories, with Lipschitz functors and Lipschitz natural transformations.

**3.4. The fundamental weighted category.** Let us come back to a δ-metric space X, and construct its fundamental weighted category, working with extended paths modulo extended 2-homotopy.

More precisely, an extended *double path* in X is a $w_\infty$-map $A: \delta \mathbf{I}^2 \to X$. It is an instance of a double homotopy (2.5), defined on the point, and the previous results for double homotopies apply; its four faces are paths in X, between four vertices. A *2-path* is a double path whose faces $\partial_1^\alpha$ are degenerate; it is a 2-homotopy $A: a \prec_2 b: x \to x'$ between its faces $\partial_2^\alpha$, which have the same endpoints. A 2-homotopy class of paths [a] is a class of the equivalence relation $\simeq_2$ spanned by the preorder $\prec_2$.

The *fundamental weighted category* $w\Pi_1(X)$ of the δ-metric space X has for objects the points of X; for arrows $[a]: x \to x'$ the 2-homotopy classes of paths from x to x', as defined above. Composition - written additively - is induced by concatenation of consecutive paths, and identities come from degenerate paths

(1)  $[a] + [b] = [a + b]$,  $\qquad 0_x = [e(x)] = [0_x]$.

The weight is defined as usual in a quotient, starting from the length of extended paths (1.8)

(2)  $w(\xi) = \inf_{a \in \xi} L(a)$.

This is indeed a category, with the same proof as for d-spaces ([G3], Thm. 3.2a), based on the constructions of 2-dimensional homotopies in 2.5a-c. Moreover, $w\Pi_1(X)$ is (sub)additively weighted, as it follows immediately from the properties of L (1.8)

(3)  $w(0_x) = 0$,  $\qquad w([a] + [b]) \le w[a] + w[b]$.

On a $\delta_\infty$-map f: X → Y, we get a $w_\infty$-functor $w\Pi_1(f)$: $w\Pi_1(X) \to w\Pi_1(Y)$

(4)  $w\Pi_1(f)(x) = f(x)$,  $\qquad w\Pi_1(f)[a] = f_*[a] = [fa]$,

$w(f_*[a]) \le \|f\|.w[a]$,  $\qquad \|f_*\| \le \|f\|$.

All this forms a functor $w\Pi_1$: $\delta_\infty$**Mtr** → $w_\infty$**Cat**, with values in the category of (small) additively weighted categories and Lipschitz functors; a functor which restricts to $\delta$**Mtr** → w**Cat**, because of the last inequality above. In particular, $w\Pi_1$ preserves Lipschitz isomorphisms and isometric isomorphisms.

Finally, a $\delta_\infty$-homotopy $\varphi$: f → g: X → Y yields a Lipschitz natural transformation

(5)  $\varphi_*$: $f_* \to g_*$: $w\Pi_1(X) \to w\Pi_1(Y)$,

$w(\varphi_*(x)) = w[\varphi(x)] \le \|\varphi\|$,  $\qquad \|\varphi_*\| \le \|\varphi\|$,

so that $w\Pi_1$: $\delta_\infty$**Mtr** → $w_\infty$**Cat** is actually a morphism of h-categories, as well as its restriction $\delta_1$**Mtr** → $w_1$**Cat** to the 1-Lipschitz case.

The fundamental weighted category of X is linked to the fundamental groupoid of the underlying space UX, by the obvious *comparison* functor

(6)  $w\Pi_1(X) \to \Pi_1(UX)$,  $\qquad x \mapsto x$,  $\quad [a] \mapsto [a]$.

**3.5. Geodesics.** In a $\delta$-metric space X, we say that an extended path a: x → x' is a *homotopic geodesic* if it realises the weight of its class, $L(a) = w[a]$, which amounts to saying that $L(a) \le L(a')$ for all extended paths $a' \simeq_2 a$.

We say that X is *geodetically simple* if every arrow $\xi$: x → x' of its fundamental weighted category $w\Pi_1(X)$ has *some* representative a which realises its weight: $L(a) = w(\xi)$; the path a is then a homotopic geodesic.

We say that X is *1-simple* if its fundamental category $w\Pi_1(X)$ is a preorder: all hom sets have at most one arrow.

The $\delta$-metric spaces $\delta\mathbf{R}^n$, $\delta\mathbf{R}^{\otimes n}$ are geodetically simple and 1-simple; all their convex subspaces are also (cf. [G3], 3.4a). The pierced plane $(!\delta\mathbf{R})^2 \setminus \{0\}$ is not geodetically simple, nor 1-simple. The $\delta$-metric sphere $\delta\mathbf{S}^1$ is geodetically simple and not 1-simple; the higher spheres are 1-simple and not geodetically simple, as antipodal points have infinitely many geodesics between them.

Being geodetically simple is somehow related with completeness of the $\delta$-metric, as it appears from these examples. Actually, a non-complete space, like $\delta]0, 1[ \subset \delta\mathbf{R}$ can be geodetically simple; but it is also true that all its extended paths  x → x'  (between two given points) stay in the compact subspace  δ[x, x'].

**3.6. Pasting Theorem** ('Seifert - van Kampen' for fundamental weighted categories). Let  X  be a δ-metric space;  $X_1, X_2$  two subspaces and  $X_0 = X_1 \cap X_2$.

If  $X = \text{int}(X_1) \cup \text{int}(X_2)$,  the following diagram of weighted categories and contracting functors (induced by inclusions) is a pushout in  w**Cat**

$$
(1) \quad
\begin{array}{ccc}
w\Pi_1 X_0 & \xrightarrow{u_1} & w\Pi_1 X_1 \\
{\scriptstyle u_2} \downarrow & \searrow{\scriptstyle v_2} & \downarrow {\scriptstyle v_1} \\
w\Pi_1 X_2 & \longrightarrow & w\Pi_1 X
\end{array}
$$

**Proof.** As in [G3], Thm.3.6. □

**3.7. Homotopy monoids.** The *fundamental weighted monoid*  $w\pi_1(X, x)$  of the δ-metric space  X  at the point  x  is the (additively) weighted monoid of endoarrows  x → x  in  $w\Pi_1(X)$.  It forms a functor from the (obvious) category  δ**Mtr**$_*$  of *pointed δ-metric spaces*, to the category of additively weighted monoids ([G8], 2.1)

(1)     $w\pi_1$: δ**Mtr**$_*$ → w$^+$**Mon**,             $w\pi_1(X, x) = w\Pi_1(X)(x, x)$.

This functor is *strictly* homotopy invariant: a *pointed homotopy*  φ: f → g: (X, x) → (Y, y)  has, by definition, a trivial path at the base-point  ($\varphi(x) = 0_y$),  whence the naturality square of every endomap  a: x → x  of  X  gives  $f_*[a] = g_*[a]$

$$
(2) \quad
\begin{array}{ccc}
f(x) & \xrightarrow{f_*[a]} & f(x) \\
{\scriptstyle [\varphi(x)]} \downarrow & & \downarrow {\scriptstyle [\varphi(x)]} \\
g(x) & \xrightarrow{g_*[a]} & g(x)
\end{array}
$$

## 4. Minimal models

This is a brief exposition of how the minimal models developed in [G7] for the fundamental category of a d-space can be enriched, in the present weighted setting. Some knowledge of the main definitions and results of [G7] would be useful.

**4.1. The fundamental weighted category of a square annulus.** Let us begin with an elementary example, enriching with a δ-metric the 'square annulus' examined in the Introduction of [G7], as an ordered topological space.

Let us start from the δ-metric space  δ**I**$^{\otimes 2}$,  with

(1)    $\delta(\mathbf{x}, \mathbf{y}) = (x_1 - y_1) + (x_2 - y_2)$    if  $x_1 \leq y_1$, $x_2 \leq y_2$,

         $= \infty$,  otherwise,

whose underlying ordered topological space (1.9.1) is the ordered topological square $\uparrow[0, 1]^2$, with euclidean topology and product order.

Taking out the *open* square $]1/3, 2/3[^2$ (marked with a cross), we get the square annulus $X \subset \delta\mathbf{I}^{\otimes 2}$, with the induced $\delta$-metric

(1) 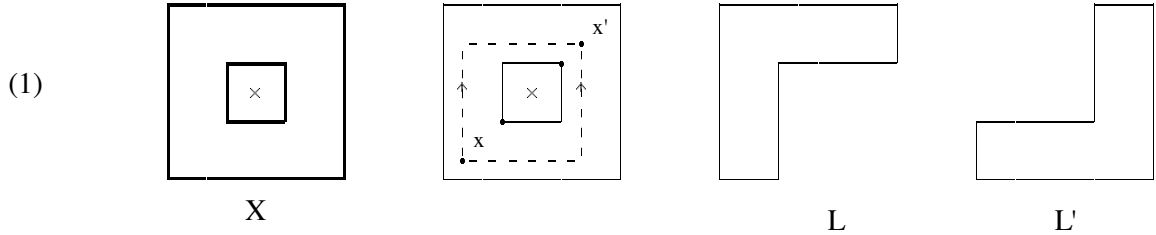

              X                                 L            L'

Its extended paths are thus the Lipschitz *order-preserving* maps $\delta[0, 1] \to X$ defined on the standard $\delta$-interval, and move 'rightward and upward' (in the weak sense). Extended homotopies of such paths are Lipschitz order-preserving maps $\uparrow[0, 1]^2 \to X$.

As a consequence of the 'van Kampen' theorem recalled above (using the subspaces L, L'), the fundamental weighted category $C = w\Pi_1(X)$ is the same as for the underlying ordered topological space: it has *some* arrow $x \to x'$ provided that $x \leq x'$ and both points are in L or L' (the closed subspaces represented above). Precisely, there are *two* arrows when $x \leq p = (1/3, 1/3)$ and $x' \geq q = (2/3, 2/3)$ (as in the second figure above), and *one* otherwise. The weight of an arrow can always be realised as the length of some representative: X is geodetically simple (3.5).

Thus, the weighted category C is 'essentially represented' by the full weighted subcategory E on four vertices 0, p, q, 1 (the central cell does not commute), *where each of the four generating arrows has weight* 2/3, and the weight of E is linear (strictly additive on composition)

(2) 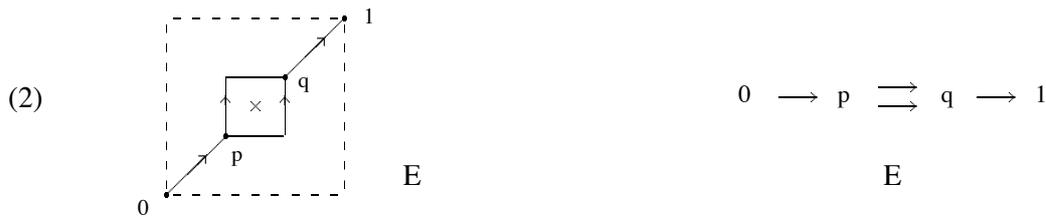

                      E                                 E

The situation can be analysed as follows, in E:

- the action begins at 0, from where we move to the point p, with weight 2/3,

- p is an (effective) future branching point, where we have to choose between two paths, each of them of weight 2/3, which join at q, an (effective) past branching point,

- from where we can only move to 1, again with weight 2/3, where the process ends.

(Definitions and properties of *regular* and *branching* points can be found in [G7], 5.3).

Recall that E is *not* equivalent to C, as a category, since C *is already a skeleton*, in the ordinary sense. In order to make precise how E can 'model' the category C, we proved in [G7] (and will recall below) that E is both *future equivalent* and *past equivalent* to C, and actually it is the 'join' of a minimal 'future model' with a minimal 'past model' of the latter. All this can now be enriched with weights.

**4.2. Future equivalence of weighted categories.** The notion of future equivalence - a symmetric version of an adjunction, *with two units* - can be easily transferred from **Cat** ([G7], 2.1) to the 2-category $w_\infty$**Cat**, since it makes sense in any 2-category.

Thus, a *future equivalence* (f, g; φ, ψ) between the weighted categories C, D will consists of a pair of Lipschitz functors and a pair of Lipschitz natural transformations, the *unit*s, satisfying two coherence conditions:

(1)   f: C $\rightleftarrows$ D :g      φ: $1_C \to gf$,   ψ: $1_D \to fg$,

(2)   fφ = ψf: f → fgf,      φg = gψ: g → gfg      (*coherence*),

and will be said to be *elementary*, or 1-Lipschitz, if both functors and both natural transformations are.

Future equivalences compose (as in [G7]), and yield an equivalence relation of categories; the elementary ones do not. Dually, *past equivalences* have *co*units, in the opposite direction.

In particular, an *elementary future retract* i: $C_0 \subset C$ will be a full weighted subcategory having a reflector p ⊣ i which is 1-Lipschitz, has a 1-Lipschitz unit η: $1_C \to$ ip and a trivial counit pi = 1. The coherence conditions of the adjunction (ηi = $1_i$, pη = $1_p$) show that the fourtuple (i, p; 1, η) is an elementary future equivalence.

A (weighted) *pf-presentation* (extending [G7], 4.2) of the weighted category C will be a diagram consisting of an elementary past retract P and an elementary future retract F of C (which are thus a full coreflective and a full reflective weighted subcategory, respectively) with elementary adjunctions $i^-$ ⊣ $p^-$ and $p^+$ ⊣ $i^+$

(3)   P $\underset{p^-}{\overset{i^-}{\rightleftarrows}}$ C $\underset{i^+}{\overset{p^+}{\rightleftarrows}}$ F      ε: $i^-p^- \to 1_C$      ($p^-i^- = 1$, $p^-ε = 1$, $εi^- = 1$),

      η: $1_C \to i^+p^+$      ($p^+i^+ = 1$, $p^+η = 1$, $ηi^+ = 1$).

**4.3. Spectra.** Coming back to the square annulus X (4.1), the weighted category C = $w\Pi_1(X)$ has a least *full reflective* weighted subcategory F, which is future equivalent to C and minimal as such. Its objects form the *future spectrum* $sp^+(C) = \{p, 1\}$ [G7]; also the full weighted subcategory F = $Sp^+(C)$ on these objects is called a *future* (weighted) *spectrum* of C

(1)

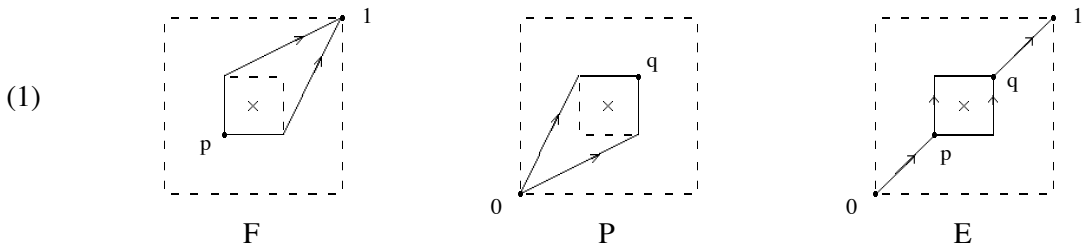

F      P      E

Dually, we have the least *full coreflective* weighted subcategory P = $Sp^-(C)$, on the *past spectrum* $sp^-(C) = \{0, q\}$.

Together, they form a (weighted) pf-presentation of C (4.2.3), called the *spectral* pf-presentation. Moreover the (weighted) *pf-spectrum* E = Sp(C) is the full weighted subcategory of C on the set of objects sp(C) = sp⁻(C)∪sp⁺(C) ([G7], 7.6). E is a *strongly minimal injective model* of the category C ([G7], Thm. 8.4).

## 5. Spaces with weighted paths

We introduce a second setting for weighted algebraic topology, which is more complicated than δ-metric spaces but has finer quotients, as we shall see in the last section. The links between the two settings are developed in Section 6.

**5.1. Main definitions.** A *w-space* X, or *space with weighted paths*, will be a topological space together with a *weight function* w: $X^I$ → [0, ∞], or *cost function* (also written $w_X$) defined on the *set* of its (continuous) paths, and satisfying three axioms concerning the constant paths $0_x$, the path-concatenation a+b of consecutive paths and *strictly* increasing reparametrisation

(wsp.0)  $w(0_x) = 0$                                                      (for all points x of X),

(wsp.1)  $w(a+b) \leq w(a) + w(b)$                        (for all pairs of consecutive paths a, b),

(wsp.2)  $w(a\rho) \leq w(a)$         (for all paths a and all *strictly increasing continuous maps* ρ: **I** → **I**).

It is easy to see that the last condition, in the presence of the others, is equivalent to asking w(aρ) ≤ w(a), for all paths a and all *increasing continuous maps* ρ: **I** → **I** which are constant on a *finite* number of subintervals. It is also equivalent to the conjunction of the following two conditions:

(wsp.1')  $w(a) \vee w(b) \leq w(a+b)$                        (for all pairs of consecutive paths a, b),

(wsp.2')  $w(a\rho) = w(a)$            (for all paths a and all *increasing homeomorphisms* ρ: **I** → **I**).

We shall say that a path is *free*, *feasible* or *unfeasible* when, respectively, its cost is 0, finite or ∞. A w-space will be said to be *linear*, or *strictly additive*, if w(a+b) = w(a) + w(b); see 5.4 for examples. We shall see in Section 6 that linear w-spaces form a coreflective subcategory. Note that we are not asking that the weight function be continuous with respect to the compact-open topology of $X^I$; in most examples, this will only be true if we restrict w to the feasible paths (or - equivalently - if we topologise [0, ∞] letting ∞ be everywhere dense, which might be interesting).

If X, Y are w-spaces, a *w-map* f: X → Y, or *map of w-spaces*, or *1-Lipschitz map* will be a continuous mapping which decreases costs: w(f∘a) ≤ w(a), for all (continuous) paths a of X. More generally, a *Lipschitz map*, or $w_\infty$-*map* f: X → Y is a continuous mapping which admits a finite Lipschitz constant λ ∈ [0, ∞[, in the sense that w(f∘a) ≤ λ.w(a), for all continuous paths a in X.

We have thus the category w**Top** of w-spaces and w-maps, embedded in the category $w_\infty$**Top** of w-spaces and $w_\infty$-maps. Again, we distinguish between *isometric isomorphisms* (of w**Top**) and *Lipschitz isomorphisms* (of $w_\infty$**Top**).

The forgetful functor U: w**Top** → **Top** has left and right adjoints D ⊣ U ⊣ C, where the discrete weight of DX is the highest possible one, with w(a) = 0 on the constant paths and ∞ on all the others, while the *natural*, codiscrete weight of CX is the lowest possible one, where all paths have a null cost. Except if otherwise stated, *when viewing a topological space as a weighted one we will use* *the embedding* C: **Top** → w**Top**, where all paths are free; this embedding is reflective (with reflector U) and coreflective (with coreflector the functor w**Top** → **Top** which sends a w-space X to its 0-ball subspace), hence preserves all limits and colimits (cf. 5.3).

Note now that the weight function of the w-space X acts on the continuous mappings a: **I** → UX, *with values in the underlying topological space*; we shall go on writing down, pedantically, such occurrences of U. (Viewing these paths as w-maps D**I** → X would also be correct, but confusing.)

Also here, reversing paths by the involution r: **I** → **I**, r(t) = 1 – t, gives the *reflected*, or *opposite*, w-space, forming a (covariant) involutive endofunctor

(1)   R: w**Top** → w**Top**,              R(X) = $X^{op}$,      $w^{op}(a)$ = w(ar).

A w-space is *symmetric* if it is invariant under reflection. It is *reflexive*, or self-dual, if it is isometrically isomorphic to its reflection (cf. 5.4). The notation X ≤ X' will mean that these w-spaces have the same underlying set and $w_X \leq w_{X'}$; equivalently, the identity of the underlying set is a w-map X' → X.

**5.2. The weight of a map.** A continuous mapping f: UX → UY between w-spaces takes paths of UX into paths of UY, and inherits thus two weights from the category of weighted sets, the *additive weight* ([G8], 1.5.3):

(1)   $|f|_0$ = $\sup_a$ (w(f∘a) – w(a))

       = min {λ | for all a, w(f∘a) ≤ λ + w(a)},              (a ∈ **Top**(**I**, UX), λ ∈ [0, ∞]).

and the *multiplicative weight*, or *Lipschitz weight* ([G8], 1.6.3):

(2)   ‖f‖ = $|f|_1$ = $\sup_a$ (w(f∘a) / w(a)) = min {λ | for all paths a, w(f∘a) ≤ λ.w(a)}.

The first distinguishes w-maps with the condition $|f|_0$ = 0. But here we shall only use the second, written as ‖f‖, which distinguishes w-maps with the condition ‖f‖ ≤ 1, and $w_\infty$-maps with the condition ‖f‖ < ∞. With this weight, w**Top** and $w_\infty$**Top** are multiplicatively weighted categories [G8]: all identities have ‖$1_X$‖ ≤ 1 and composition gives ‖gf‖ ≤ ‖f‖.‖g‖.

If X is a w-space and λ ∈ [0, ∞[, we write λX the same topological space equipped with the weight λ.$w_X$. A $w_\infty$-map f: X → Y with ‖f‖ ≤ λ is the same as a w-map λX → Y.

**5.3. Limits.** The category w**Top** has all limits and colimits, computed as in **Top** and equipped with a suitable w-structure.

Thus, for a product Π$X_i$, a path a: **I** → U(Π$X_i$) of components $a_i$: **I** → U$X_i$ has weight w(a) = sup w($a_i$). For a sum Σ$X_i$, a path a: **I** → U(Σ$X_i$) lives in *one* component U$X_i$ and inherits the weight from the latter.

Given a pair of parallel w-maps f, g: X → Y, the equaliser is the topological one, with the restricted weight function. The coequaliser is the topological coequaliser Y/R, with the induced weight

(1)   w(c) = inf($\sum_i w_Y(b_i)$)                                            (c: **I** → U(Y/R)),

the inf being taken on all finite families ($b_1$,..., $b_n$) of paths in UY such that their projections on the quotient, p$b_i$: **I** → U(Y/R), are consecutive, and give c = ((p$b_1$) + ... + (p$b_n$))ρ  (by n-ary concatenation and reparametrisation by an increasing homeomorphism $\mathbf{I} \to \mathbf{I}$). Of course, if there are no such families, $w(c) = \inf(\emptyset) = \infty$.

The only non-trivial points are verifying that this weight-function satisfies the axioms (wsp.1, 1') of 5.1. Take $c = c' + c'': \mathbf{I} \to Y/R$. First, every pair of decompositions of $c'$, $c''$

(2)  $c' = ((pb_1') + \ldots + (pb_m'))\rho'$,  $\qquad c'' = ((pb_1'') + \ldots + (pb_n''))\rho''$,

gives a decomposition $c = ((pb_1') + \ldots + (pb_n''))\rho$; therefore

(3)  $(w_Y(b_1') + \ldots + w_Y(b_m')) + (w_Y(b_1'') + \ldots + w_Y(pb_n'')) \geq w(c)$,

and $w(c') + w(c'') \geq w(c)$. Second, given a decomposition $c = ((pb_1) + \ldots + (pb_n))\rho$, we can always assume that $n = 2k$ (possibly inserting a constant path, without modifying $\Sigma_i w_Y(b_i)$). Then

(4)  $((pb_1) + \ldots + (pb_k)) + ((pb_{k+1}) + \ldots + (pb_{2k})) = c\rho^{-1} = (c'\rho') + (c''\rho'')$,

   $c'\rho' = (pb_1) + \ldots + (pb_k)$,

   $w(c') \leq w_Y(b_1) + \ldots + w_Y(b_k) \leq w_Y(b_1) + \ldots + w_Y(b_{2k})$,

(for suitable reparametrisations $\rho'$, $\rho''$). It follows that $w(c') \leq w(c)$, and similarly for $c''$.

Linear w-spaces are not closed under (even binary) products, as we see below. But they are closed under subspaces (plainly), all colimits (by adjointness) and tensor product (5.5).

The category $w_\infty\mathbf{Top}$ has finite limits and colimits, which can be constructed as above. Even if, now, they are only determined up to *Lipschitz* isomorphism, we shall keep the previous construction as privileged ones. Thus, when we write $X \times Y$ in $w_\infty\mathbf{Top}$, we still mean that its weight is the $l_1$-weight, with $w(a, b) = w(a) + w(b)$; isomorphic constructions will have different names (cf. 5.5). It is easy to verify that: $\|f \times g\| \leq \|f\| \vee \|g\|$.

**5.4. Standard models.** The *standard weighted real line*, or *w-line* $w\mathbf{R}$, is the euclidean line with the following weight on all paths $a: \mathbf{I} \to \mathbf{R}$, equivalently defined by its span or length in $\delta\mathbf{R}$ (1.8)

(1)  $w(a) = sp(a) = L(a)$;

thus, $w(a)$ is finite if and only if $a$ is a (weakly) *increasing* path, and then $w(a) = a(1) - a(0)$. (General links between $\delta$-metric spaces and w-spaces will be studied in the next section.)

Its cartesian power in $w\mathbf{Top}$, the n-*dimensional real w-space* $w\mathbf{R}^n$ has $w(a) = \sup_i(a_i(1) - a_i(0))$ for all increasing paths $a: \mathbf{I} \to \mathbf{R}^n$ (with respect to the product order, $x \leq x'$ if and only if $x_i \leq x_i'$ for all i). Plainly, $w\mathbf{R}$ is linear while every higher dimensional $w\mathbf{R}^n$ is not.

The *standard w-interval* $w\mathbf{I}$ has the subspace structure of the w-line; the *standard w-cube* $w\mathbf{I}^n$ is its n-th power, and a subspace of $w\mathbf{R}^n$. These w-spaces are not symmetric (for $n > 0$), but reflexive; in particular, the canonical reflecting isomorphism

(2)  $r: w\mathbf{I} \to (w\mathbf{I})^{op}$, $\qquad t \mapsto 1 - t$,

will be used to *reflect* paths and homotopies.

The *standard weighted circle* $w\mathbf{S}^1$ will be the coequaliser in $w\mathbf{Top}$ of the following two pairs of maps (equivalently)

(3)  $\partial^-, \partial^+: \{*\} \rightrightarrows w\mathbf{I}$, $\qquad \partial^-(*) = 0, \quad \partial^+(*) = 1$,

(4)   id, f: w**R** $\rightrightarrows$ w**R**,                              f(x) = x + 1.

Thus, the 'standard realisation' of the first coequaliser is the quotient (w**I**)/∂**I**, which identifies the endpoints; a feasible path turns around the circle in a precise direction, and its weight measures the length of the path with respect to the length of the circle: w(a) = L(a) in δ**S**$^1$. The Lipschitz-isomorphic structure 2π.w**S**$^1$ is also of interest. Both are linear.

More generally, the *weighted n-dimensional sphere* will be the quotient of the weighted cube w**I**$^n$ modulo its (ordinary) boundary ∂**I**$^n$, while w**S**$^0$ has the discrete topology and the unique w-structure

(5)   w**S**$^n$ = (w**I**$^n$)/(∂**I**$^n$)   (n > 0),                  w**S**$^0$ = **S**$^0$ = {−1, 1}.

All directed spheres are reflexive. Again, w**S**$^1$ is linear while the higher spheres are not.

**5.5. Tensor product.** The tensor product X⊗Y of two w-spaces (similar to the tensor product in **w**$^+$) will be the cartesian product of the underlying topological spaces, with an $l_1$-type weight (instead of the $l_\infty$-one, pertaining to the cartesian product)

(1)   $w_\otimes(a, b) = w_X(a) + w_Y(b)$,

where (a, b): **I** → X×Y denotes the path of components a: **I** → X, b: **I** → Y. Plainly, this defines a symmetric monoidal structure on w**Top**, with identity the singleton space {∗}.

Linear w-spaces are closed under tensor product. In particular, all tensor powers (w**R**)$^{\otimes n}$, (w**I**)$^{\otimes n}$ and (w**S**$^1$)$^{\otimes n}$ are linear. The following theorem shows that all of them are exponentiable with respect to the tensor product; in particular, the tensor power (w**I**)$^{\otimes n}$, which is what is relevant for homotopy.

This tensor product extends to $w_\infty$**Top**, with ‖f⊗g‖ ≤ ‖f‖v‖g‖. Even if here the tensor product is isomorphic to the cartesian one, we will keep these two construct distinct.

**5.6. Theorem** (Exponentiable w-spaces). Let Y be a *linear w-space with a locally compact Hausdorff topology*. Then Y is ⊗-*exponentiable* in w**Top**; for every w-space Z, the internal hom

(1)   $Z^Y$ = w**Top**(Y, Z) ⊂ **Top**(UY, UZ),

is the set of w-maps, with the compact-open topology restricted from (UZ)$^{UY}$ and the w-structure where a path c: **I** → U($Z^Y$) ⊂ (UZ)$^{UY}$ has the following weight

(2)   $W(c) = \sup_b (w_Z(ev \circ (c, b)) - w_Y(b))$,                            (for b: **I** → UY),

where λ − μ is the truncated difference in **w**$^+$, and the evaluation mapping

(3)   ev: $Z^Y$ ⊗ Y → Z,

is the restriction of the topological one; this mapping is a w-map, the counit of the adjunction.

**Proof.** (Note. The same proof, conveniently simplified, shows that δ**Mtr** is monoidal closed.)

We defer to the end the technical part showing that (1) is indeed a w-structure. First, the evaluation mapping (3) satisfies the inequality $w_\otimes(c, b) \geq w_Z(ev \circ (c, b))$, because of the adjunction in **w**$^+$:

(4)   $W(c) \geq w_Z(ev \circ (c, b)) - w_Y(b)$                          (for all b: **I** → UY),

   $w_\otimes(c, b) = W(c) + w_Y(b) \geq w_Z(ev \circ (c, b))$.

Second, the pair $(Z^Y, \text{ev}: Z^Y \otimes Y \to Z)$ is a universal arrow from the functor $-\otimes Y$ to the object Z: given a w-space X and a w-map $f: X \otimes Y \to Z$, we have to prove that there is precisely one w-map $g: X \to Z^Y$ such that f factors as:

(5) $\quad \text{ev} \circ (g \otimes Y): X \otimes Y \to Z^Y \otimes Y \to Z.$

Indeed, since Y is exponentiable in **Top**, there exists precisely one continuous mapping $g: UX \to (UZ)^{UY}$ such that $f = \text{ev} \circ (g \times UY)$, and it will be sufficient to prove the following two facts.

(a) $\text{Im}(g) \subset Z^Y$. For $x \in X$, we must prove that $g(x): Y \to Z$ is a w-map. And indeed, for every path $b: \mathbf{I} \to UY$

(6) $\quad w(g(x) \circ b) = w(f \circ (0_x, b)) \leq w_\otimes(0_x, b) = w_X(0_x) + w_Y(b) = w_Y(b).$

(b) The mapping g is a w-map $X \to Z^Y$. And indeed, for every path $a: \mathbf{I} \to UA$

(7) $\quad W(g \circ a) = \sup_b (w_Z(\text{ev} \circ (ga, b)) - w_Y(b)) = \sup_b (w_Z(f \circ (a, b)) - w_Y(b))$

$\quad \leq \sup_b (w_\otimes(a, b) - w_Y(b)) = w_X(a).$

Finally, we verify the axioms for the hom-weight W. First, the constant path $0_h: \mathbf{I} \to Z^Y$ at some w-map $h: Y \to Z$ gives

(8) $\quad W(0_h) = \sup_b (w_Z(\text{ev} \circ (0_h, b)) - w_Y(b)) = \sup_b (w_Z(hb) - w_Y(b)) = 0.$

Second, to prove (wsp.1), let $c = c' + c''$ be a concatenation of paths in $U(Z^Y)$. We can always rewrite a path $b: \mathbf{I} \to UY$ as the concatenation $b = b' + b''$ of its two halves, so that, using the assumption that Y *is linear*:

(9) $\quad W(c' + c'') = \sup_b ((w_Z(\text{ev} \circ (c' + c'', b)) - w_Y(b))$

$\quad = \sup_{b'b''} (w_Z(\text{ev} \circ (c' + c'', b' + b'')) - w_Y(b' + b''))$ (for all consecutive paths b', b'' in UY),

$\quad = \sup_{b'b''} (w_Z((\text{ev} \circ (c', b')) + (\text{ev} \circ (c'', b''))) - w_Y(b') - w_Y(b''))$

$\quad \leq \sup_{b'b''} (w_Z(\text{ev} \circ (c', b')) - w_Y(b')] + [w_Z(\text{ev} \circ (c'', b'')) - w_Y(b'')] \leq W(c') + W(c''),$

since the last term amounts to the previous sup for arbitrary paths b', b'' in Y. (Note. For δ-metric spaces, one would use the fact that all 'paths' $(y, y'): \mathbf{2} \to Y$ *can* be rewritten as a trivial 'concatenation' $(y, y') + (y', y')$, with $d(y, y') = d(y, y') + d(y', y')$.)

Now, for (wsp.1'), we can make our least upper bound smaller by restriction to those paths $b: \mathbf{I} \to Y$ which are constant on $[1/2, 1]$, so that $b = b' + b''$ with an arbitrary b' and b'' constant at the terminal of b':

(10) $\quad W(c' + c'') \geq \sup_{b'} (w_Z((\text{ev} \circ (c', b')) + (\text{ev} \circ (c'', b''))) - w_Y(b' + b''))$

$\quad \geq \sup_{b'} (w_Z(\text{ev} \circ (c', b')) - w_Y(b')) = W(c').$

where, again, we have used the linear property of Y: $w(b) = w(b') + w(b'') = w(b')$.

Last, for (wsp.2'), given an increasing homeomorphism $\rho: \mathbf{I} \to \mathbf{I}$, every path b' in UY can be rewritten as $b\rho$, with $b = b'\rho^{-1}$, so that:

(11) $\quad W(c\rho) = \sup_b (w_Z(\text{ev} \circ (c\rho, b\rho)) - w_Y(b\rho)) = \sup_b (w_Z(\text{ev} \circ (c, b) \circ \rho) - w_Y(b\rho)) = W(c). \quad \square$

**5.7. Elementary and extended paths.** Let $X$ be a w-space. An *elementary path* (resp. an *extended path*, or *Lipschitz path*) in $X$ will be a 1-Lipschitz (resp. a Lipschitz) map $a: w\mathbf{I} \to X$.

Thus, a continuous mapping $a: \mathbf{I} \to X$ is an elementary path if and only if $\|a\| \leq 1$, for the Lipschitz weight (5.2.2)

(1)   $\|a\| = \min\{\lambda \in [0, \infty] \mid \text{for all increasing maps } \rho: \mathbf{I} \to \mathbf{I},\ w(a\rho) \leq \lambda.(\rho(1) - \rho(0))\}.$

and is an extended path if and only if $\|a\| < \infty$. Elementary paths cannot be concatenated, because - loosely speaking - this procedure doubles the velocity, whose least upper bound is the Lipschitz weight.

The *reflected* (elementary or extended) path is obtained in the obvious way

(2)   $a^{op} = ar: w\mathbf{I} \to X^{op},$       $r(t) = 1 - t,$

A *reversible* extended path is a mapping $a: \mathbf{I} \to X$ such that both $a$ and $a^{op}$ are extended paths $w\mathbf{I} \to X$.

In the category w**Top**, pasting two copies of the standard weighted interval, one after the other, can be realised as $w[0, 2] \subset w\mathbf{R}$ (or as $2.w\mathbf{I}$, cf. 5.2), which is of no help to concatenate paths parametrised on $w\mathbf{I}$, in w**Top**. But in $w_\infty$**Top** this pasting can be realised as $w\mathbf{I}$ (Lipschitz-isomorphic to $w[0, 2]$), by the *standard concatenation pushout*

(3)   $\begin{array}{ccc} \{*\} & \xrightarrow{\partial^+} & w\mathbf{I} \\ \partial^- \downarrow & & \downarrow k^- \\ w\mathbf{I} & \xrightarrow[k^+]{} & w\mathbf{I} \end{array}$       $k^-(t) = t/2, \quad k^+(t) = (t + 1)/2.$

(Note that the diagram above is still in w**Top**, but is a pushout only in $w_\infty$**Top**.) Now, given two consecutive Lipschitz paths $a, b: w\mathbf{I} \to X$, with $a(1) = b(0)$, we get a concatenated path

(4)   $a + b: w\mathbf{I} \to X,$       $\|a + b\| \leq 2.(\|a\| + \|b\|),$

as follows from the following proposition (or using the pushout $2.w\mathbf{I}$, in w**Top**).

We can now treat homotopies as for δ-metric spaces, in Section 2, distinguishing the h-categories $w_1$**Top** $\subset$ w**Top** $\subset w_\infty$**Top**, where w**Top** has maps in w**Top** and homotopies in $w_\infty$**Top**. And define the fundamental weighted category of a w-space as in Section 3. In particular, concatenation is based on the following result.

**5.8. Proposition.** For every w-space $X$, the functor $X \times -: w_\infty\mathbf{Top} \to w_\infty\mathbf{Top}$ preserves the standard concatenation pushout (5.7.3). Moreover, if a map $f: X \times w\mathbf{I} \to Y$ comes from the pasting of two 'consecutive' maps $f_0, f_1: X \times w\mathbf{I} \to Y$, we have the upper bound for its Lipschitz weight:

(1)   $\|f\| \leq 2.(\|f_0\| + \|f_1\|)$       $(f_0 = f \circ (X \times k^-),\ f_1 = f \circ (X \times k^+)).$

Equivalently, one can use the Lipschitz-isomorphic functor $X \otimes -$.

**Proof.** In **Top**, the preservation holds because the subspaces $UX \times [0, 1/2]$ and $UX \times [1/2, 1]$ form a finite closed covering of $UX \times \mathbf{I}$, so that each mapping defined on the latter and continuous on such closed parts is continuous.

Consider then a (topological) map $f: UX \times I \to UY$ coming from the pasting of two maps $f_0, f_1$ on the topological pushout $UX \times I$

(2) $f(x, t) = f_0(x, 2t)$, for $0 \leq t \leq 1/2$, $\quad\quad\quad f(x, t) = f_1(x, 2t - 1)$, for $1/2 \leq t \leq 1$.

Let now $(a, \rho): I \to UX \times I$ be any feasible path; in particular, $\rho: I \to I$ is an increasing map. If the image of $\rho$ is contained in the first half of $I$, then $f \circ (a, \rho) = f_0(a, 2\rho)$ and

(3) $w(f \circ (a, \rho)) \leq \|f_0\|.(w(a) \vee 2w(\rho)) \leq 2.\|f_0\|.w(a, \rho)$.

Similarly for the second half. Otherwise, since $\rho$ is increasing, we have $\rho(t_1) = 1/2$ at some interior point $t_1 \in\ ]0, 1[$; and we can *assume* that $t_1 = 1/2$ (up to precomposing with an *increasing homeomorphism* $\sigma: I \to I$, which does not modify the weight of paths, by (wsp.2')). Now, the path $f \circ (a, \rho): I \to UY$ is the concatenation of two paths $c_i: wI \to UY$ which factor through the Lipschitz maps $f_i$

(4) $c_0(t) = f \circ (a(t/2), \rho(t/2)) = f_0 \circ (a(t/2), 2\rho(t/2))$,

$c_1(t) = f \circ (a((t+1)/2), \rho((t+1)/2)) = f_1 \circ (a((t+1)/2), 2\rho((t+1)/2) - 1)$.

and finally we can conclude that $f$ is Lipschitz, with the upper bound (1)

(5) $w(f \circ (a, \rho)) \leq w(c_0) + w(c_1) \leq (\|f_0\| + \|f_1\|).(w(a) \vee 2w(\rho)) \leq 2.(\|f_0\| + \|f_1\|).w(a, \rho)$. $\quad\square$

## 6. Linear and metrizable w-spaces

The *span* and *length* function of a $\delta$-metric space $X$, defined in 1.8, allow us to construct the w-spaces $\mathbf{sp}X$ (6.2) and $\mathbf{L}X$; the latter is linear (6.3).

**6.1. Linear w-spaces.** First, we want to observe that linear w-spaces form a full subcategory $\mathrm{Lw}_\infty\mathbf{Top}$ of $\mathrm{w}_\infty\mathbf{Top}$, which has a *coreflector* L, right adjoint to the embedding U

(1) $U: \mathrm{Lw}_\infty\mathbf{Top} \rightleftarrows \mathrm{w}_\infty\mathbf{Top} : L$ $\quad\quad\quad\quad\quad\quad\quad\quad\quad\quad\quad (U \dashv L)$.

In fact, for a w-space $X$, there is a *linearised* w-space $L(X)$ on the same underlying topological space, endowed with the least linear weight $L \geq w$

(2) $L(a) = \sup_{\mathbf{t}} \Sigma_i\, w(a(t_{i-1}, t_i))$ $\quad\quad\quad\quad\quad\quad\quad\quad\quad (0 = t_0 < t_1 < ... < t_p = 1)$,

$a(t_{i-1}, t_i)(t) = a((1 - t).t_{i-1} + t.t_i)$ $\quad\quad\quad\quad\quad\quad\quad\quad\quad\quad\quad\quad\quad (0 \leq t \leq 1)$.

Note that we have written $a(t_{i-1}, t_i): I \to X$ the restriction of the path $a$ to the interval $[t_{i-1}, t_i]$, reparametrised on the standard interval.

Thus $L(X) \geq X$, and $L(X) = X$ if and only if the w-space $X$ is linear. These relations form, respectively, the counit $UL \to 1$ and the unit $LU = 1$ of the adjunction.

All this restricts to contractions, giving a full coreflective subcategory $\mathrm{Lw}\mathbf{Top}$ of $\mathrm{w}\mathbf{Top}$. Therefore, linear w-spaces are closed under colimits in $\mathrm{w}\mathbf{Top}$.

**6.2. Span-metrizable w-spaces.** Now, let us construct an adjunction

(1)  $\delta: \mathrm{w}_\infty\mathbf{Top} \rightleftarrows \delta_\infty\mathbf{Mtr} :\mathbf{sp}$                    ($\delta \dashv \mathbf{sp}$).

First, the functor

(2)  $\delta: \mathrm{w}_\infty\mathbf{Top} \to \delta_\infty\mathbf{Mtr}$,              $\|\delta f\| \leq \|f\|$,              ($\delta: \mathrm{w}\mathbf{Top} \to \delta\mathbf{Mtr}$),

sends a w-space $X$ to the δ-metric space $\delta X$, consisting of the same set with

(3)  $\delta(x, x') = \inf_a w(a)$,

where $a: \delta\mathbf{I} \to X$ varies in the set of extended paths in $X$, from $x$ to $x'$. The δ-metric spaces obtained in this way from a w-space will be said to be *geodetic*. Plainly, if $f: X \to X'$ is a $\mathrm{w}_\infty$-map, $\delta f = f: \delta X \to \delta X'$ is continuous and satisfies the inequality of (2), whence it is a $\delta_\infty$-map (and 1-Lipschitz if $f$ is).

Second, the functor

(4)  $\mathbf{sp}: \delta_\infty\mathbf{Mtr} \to \mathrm{w}_\infty\mathbf{Top}$,              $\|\mathbf{sp}f\| \leq \|f\|$,              ($\mathbf{sp}: \delta\mathbf{Mtr} \to \mathrm{w}\mathbf{Top}$),

has essentially been constructed in Section 1. For a δ-metric space $Y$, we let $\mathbf{sp}Y$ be the same set equipped with the symmetric topology (1.6) and with the weight-function sp (1.8.1)

(5)  $\mathrm{sp}(a) = \sup_t \delta(a(t_0), a(t_1))$                    ($0 \leq t_0 < t_1 \leq 1$),

which has already been seen to satisfy the axioms of w-spaces (1.8).

The w-spaces obtained in this way will be said to be *span-metrizable*. On maps, we take again the same underlying mapping.

These two functors form an idempotent adjunction $\delta \dashv \mathbf{sp}$, which restricts to a (covariant) Galois connections whenever we fix the underlying set: in fact, both functors do not change the latter, and unit and counit reduce to inequalities

(6)  $X \geq \mathbf{sp}(\delta X)$,                    $\delta(\mathbf{sp}Y) \geq Y$,

where $X$ is a w-space and $Y$ a δ-metric space. (For idempotent adjunctions, see [AT] Section 6 and [LS] Lemma 4.3.)

The adjunction gives an equivalence between the full subcategories of:

(a) span-metrizable w-spaces, characterised by the condition $X = \mathbf{sp}(\delta X)$, or equivalently by the condition $X = \mathbf{sp}(Y)$ for a suitable δ-metric structure $Y$ (on the same set),

(b) geodetic δ-metric spaces, characterised by the condition $Y = \delta(\mathbf{sp}Y)$, or equivalently by the condition $Y = \delta(X)$ for some weighted structure $X$ on the associated topological space.

Restricting to 1-Lipschitz maps, span-metrizable w-spaces form a reflective subcategory of w**Top**, closed under limits, while geodetic δ-metric spaces form a coreflective subcategory of δ**Mtr**, closed under colimits.

In the standard examples of 5.4, the standard w-line is span-metrizable, $\mathrm{w}\mathbf{R} = \mathbf{sp}(\delta\mathbf{R})$ and the standard δ-line is geodetic, $\delta(\mathrm{w}\mathbf{R}) = \delta\mathbf{R}$. Similarly, in higher dimension, $\mathrm{w}\mathbf{R}^n = \mathbf{sp}(\delta\mathbf{R}^n)$ and $\delta(\mathrm{w}\mathbf{R}^n) = \delta\mathbf{R}^n$; this also holds for the standard interval and its powers.

The standard δ-circle $\delta\mathbf{S}^1 = \delta(\mathrm{w}\mathbf{S}^1)$ is geodetic, while the circle $\mathbf{S}^1$ with the euclidean metric of $\mathbf{R}^2$ is not, since $\delta(\mathbf{sp}\mathbf{S}^1)$ has the obvious geodetic distance, which is bigger. The standard w-circle $w\mathbf{S}^1$ is not span-metrizable, since the weight (i.e. length) of its feasible paths has no finite upper bound, while the δ-metric of $\delta\mathbf{S}^1 = \delta(w\mathbf{S}^1)$ cannot exceed 1.

**6.3. The length adjunction.** The span-adjunction (6.2.1) and the adjunction of linear w-spaces (6.1.1) give a composed adjunction, idempotent again

(1)   $\delta: \mathbf{Lw}_\infty\mathbf{Top} \rightleftarrows \delta_\infty\mathbf{Mtr} :\mathbf{L}$                                                         $(\delta \dashv \mathbf{L})$.

Now, $\delta$ is the restriction of the functor 6.2.2, and equips a linear w-space X with the geodetic δ-metric $\delta(x, x') = \inf_a w(a)$. On the other hand, $\mathbf{L} = L \circ sp$ takes a δ-metric space Y to the same set equipped with the symmetric topology (1.6) and with the (linear) weight-function L which we have already defined in 1.8

(2)   $L(a) = \sup_\mathbf{t} \Sigma_i\, \delta(a(t_{i-1}), a(t_i))$                                        $(0 = t_0 < t_1 < ... < t_p = 1)$.

Also here, maps are left 'unchanged' and $\|\delta f\| \leq \|f\|$, $\|Lf\| \leq \|f\|$, so that the adjunction restricts to contractions.

**6.4. Length-metrizable w-spaces.** Also the length adjunction (6.3.1) becomes a (covariant) Galois connection whenever we fix the underlying set: unit and counit reduce to inequalities

(1)   $X \geq \mathbf{L}(\delta X)$,                                        $\delta(\mathbf{L}Y) \geq Y$,

where X is a linear w-space and Y a δ-metric space.

The adjunction gives thus an equivalence between the full subcategories of:

(a) *length-metrizable w-spaces*, characterised by the condition $X = \mathbf{L}(\delta X)$, or equivalently by the condition $X = \mathbf{L}Y$ for some δ-metric structure Y on the same set (all such w-spaces are linear),

(b) *linearly geodetic* δ-metric spaces, characterised by the condition $Y = \delta(\mathbf{L}Y)$, or $Y = \delta X$ for some *linear* weight X on the associated topological space.

Thus, a linearly geodetic δ-metric space is geodetic; the converse need not be true. For instance, the δ-metric subspace $Y \subset \delta\mathbf{R}^2$ consisting of the union of the axes is geodetic, but not linearly so: the points $y = (-1, 0)$ and $y' = (0, 1)$ have $\delta(y, y') = 1$ but all feasible paths a in Y, from y to y', have length $L(a) = 2$.

The two notions of 'metrizability' of w-spaces are not comparable. Indeed, the w-line $w\mathbf{R}$ is metrizable in both senses. The w-plane $w\mathbf{R}^2$ is span-metrizable and not linear, hence not length-metrizable. The standard w-circle $w\mathbf{S}^1$ (5.4) is only length-metrizable. Finally, in 7.2, we will show that the irrational rotation w-space $W_\vartheta$ has a trivial δ-metric on $\delta W_\vartheta$ (always zero), whence it is neither span- nor length-metrizable.

**6.5. Directed spaces.** Finally, we have a forgetful functor

(1)   $\mathbf{d}: w_\infty\mathbf{Top} \to d\mathbf{Top}$,

which sends a w-space to the same topological space, equipped with the feasible paths as distinguished ones. Composing $\mathbf{L}: \delta_\infty\mathbf{Mtr} \to \mathbf{Lw}_\infty\mathbf{Top}$ with the latter, we get the forgetful functor $\delta_\infty\mathbf{Mtr} \to d\mathbf{Top}$ already considered in 1.9 (which distinguishes the L-feasible paths of a δ-metric space).

## 7. Weighted noncommutative tori and their classification

The *irrational rotation w-space* $W_\vartheta$ (7.1.3) has classifications similar to the irrational rotation C*-algebras $A_\vartheta$. Analogous results have been obtained in [G6] for 'normed' cubical sets and their 'normed' homology - an earlier occurrence of weighted algebraic topology; cubical sets give weaker results, without the metric aspects (cf. [G5]). Throughout this section $\vartheta$ is an irrational number.

**7.1. Irrational rotation w-spaces.** Let us begin recalling some well-known 'noncommutative spaces'.

First, take the line $\mathbf{R}$ and its (dense) additive subgroup $G_\vartheta = \mathbf{Z} + \vartheta\mathbf{Z}$, acting on the former by translations. In **Top**, the orbit space $\mathbf{R}/G_\vartheta = S^1/\vartheta\mathbf{Z}$ is trivial: an uncountable set with the coarse topology. In $\delta$**Mtr**, the quotient $\delta\mathbf{R}/G_\vartheta$ is trivial as well: an uncountable set with the null 'metric'.

In noncommutative geometry, this set is 'interpreted' as the (noncommutative) C*-algebra $A_\vartheta$, generated by two unitary elements u, v under the relation $vu = \exp(2\pi i\vartheta).uv$, and called the *irrational rotation algebra* associated with $\vartheta$, or also a *noncommutative torus* [C1, C2, Ri, PV]. Both its complex K-theory groups are two-dimensional.

A relevant achievement of K-theory [PV, Ri] classifies these algebras, by proving that $K_0(A_\vartheta) \cong \mathbf{Z} + \vartheta\mathbf{Z}$ *as an ordered subgroup of* $\mathbf{R}$; moreover, the traces of the projections of $A_\vartheta$ cover the set $G_\vartheta \cap [0, 1]$. It follows that $A_\vartheta$ and $A_{\vartheta'}$ are *isomorphic* if and only if $\vartheta' \in \pm\vartheta + \mathbf{Z}$ [Ri, Thm. 2] and *strongly Morita equivalent* if and only if $\vartheta$ and $\vartheta'$ are equivalent modulo the *fractional* action (on the irrationals) of the group $GL(2, \mathbf{Z})$ of invertible integral 2×2 matrices [Ri, Thm. 4]

(1) $\quad \begin{pmatrix} a & b \\ c & d \end{pmatrix}.t = \dfrac{at+b}{ct+d} \qquad\qquad (a, b, c, d \in \mathbf{Z};\ ad - bc = \pm 1),$

(or the action of the projective general linear group $PGL(2, \mathbf{Z})$ on the projective line). Since the group $GL(2, \mathbf{Z})$ is generated by the matrices

(2) $\quad R = \begin{pmatrix} 0 & 1 \\ 1 & 0 \end{pmatrix}, \qquad T = \begin{pmatrix} 1 & 1 \\ 0 & 1 \end{pmatrix},$

the orbit of $\vartheta$ is its *closure* $\{\vartheta\}_{RT}$ under the transformations $R(t) = t^{-1}$ and $T^{\pm 1}(t) = t \pm 1$ (on $\mathbf{R}\setminus\mathbf{Q}$)

We show now how one can obtain similar results with the w-space naturally arising from the action of $G_\vartheta$ *on the w-line*: the point is to replace a topologically-trivial orbit space $\mathbf{R}/G$ with the corresponding quotient of the standard w-line $w\mathbf{R}$, by a procedure analogous to the one followed in [G5, G6] for cubical sets or weighted cubical sets.

We are thus lead to consider the *irrational rotation w-space*

(3) $\quad W_\vartheta = (w\mathbf{R})/G_\vartheta,$

whose feasible paths reduce to the *projection of the feasible paths* of $w\mathbf{R}$, as we prove now. It will be useful to use the following standard weight on the additive groups $\mathbf{R}$ and $G_\vartheta$

(4) $\quad w(x) = \delta(0, x),$

i.e. $w(x) = x$ when $x \geq 0$, $w(x) = \infty$ otherwise. And the (restricted) standard weight $w(x) = x$ on the additive monoids $\mathbf{R}^+$ and $G_\vartheta^+ = G_\vartheta \cap \mathbf{R}^+$ formed by the elements of finite weight.

**7.2. Theorem.** (a) The fundamental weighted monoid of $W_\vartheta$ at each point $\bar{x} \in \mathbf{R}/G_\vartheta$ is isometrically isomorphic to the additive weighted monoid $G_\vartheta^+$, via the weight function

(1) $\quad$ w: $w\pi_1(W_\vartheta, \bar{x}) \to [0, \infty[, \qquad\qquad \mathrm{Im}(w) = G_\vartheta^+.$

(b) Choosing a representative $x \in \mathbf{R}$ of $\bar{x}$, for every *feasible* path $\bar{a}: w\mathbf{I} \to W_\vartheta$ starting at $\bar{x}$ there is precisely one increasing path $a: w\mathbf{I} \to \mathbf{R}$ which lifts it and starts at $x$. Moreover, the weight of $\bar{a}$ in $W_\vartheta$ coincides with the weight of $a$, $w(a) = a(1) - a(0) = a(1) - x$.

(c) The w-space $W_\vartheta$ is linear; the associated metric space $\delta W_\vartheta$ (6.2.2) is codiscrete, with $\delta(\bar{x}, \bar{y})$ always zero, so that $W_\vartheta$ is neither span- nor length-metrizable.

**Proof.** We begin with proving (b). Take a feasible path $\bar{a}: w\mathbf{I} \to W_\vartheta$ starting at $\bar{x}$, and choose a representative $x \in \mathbf{R}$ of the latter. There exists then some finite family $a_1,\ldots, a_p$ of feasible (i.e., increasing) paths in $w\mathbf{R}$ such that the projections $\bar{a}_j$ are consecutive and give $\bar{a} = \bar{a}_1 + \ldots + \bar{a}_p$; further, $w(\bar{a})$ is the greatest lower bound of $\Sigma\, w(a_j)$ (for such families).

Now, up to $G_\vartheta$-translations, we may assume that $a_1$ starts at $x$ and all $a_j$ are consecutive (without changing their weight and the concatenation of projections). Thus $a = a_1 + \ldots + a_p$ projects to $\bar{a}$ with $w(a) = \Sigma\, w(a_j)$, since $w\mathbf{R}$ is linear. It follows that $w(\bar{a})$ is the greatest lower bound of $w(a)$, where $a$ varies among the paths in $\mathbf{R}$ which start at $x$ and lift $\bar{a}$.

But there is only one path which satisfies these conditions. Indeed, if also $b$ does, the image of the continuous mapping $a - b: \mathbf{I} \to \mathbf{R}$ must be contained in $G_\vartheta$, which is totally disconnected; thus $a - b$ is constant, and $a(0) = x = b(0)$ gives $a = b$. It follows that $w(\bar{a}) = w(a)$, where $a$ is the unique path in $\mathbf{R}$ which starts at $x$ and lifts $\bar{a}$.

For (c), we have proved above that $W_\vartheta$ is linear. The other assertions are obvious, taking into account the characterisations of span- and length-metrizable w-spaces, in 6.2, 6.4.

For (a), let us consider the weight function (1). First, we show that its image is $G_\vartheta^+$. For a loop $\bar{a}$, we have $w(\bar{a}) = w(a)$ where the (increasing) lifting $a$ starts at $x$ and ends at some $x' \geq x$, which also projects to $\bar{x}$; thus $w(\bar{a}) = x' - x \in G_\vartheta^+$. On the other hand, if $g \in G_\vartheta^+$, any increasing path $a: x \to x + g$ projects to a loop at $\bar{x}$, whose weight is $g$.

Finally, we must prove that the weight function is injective. Let $\bar{a}, \bar{b}$ be two loops at $\bar{x}$ with the same weight $g \in G_\vartheta^+$, and let $a, b$ be their lifting starting at $x$; they have again the same weight $g$, which means that they end at the same point $x' = x + g$. Then, the increasing path $c = a \vee b: \mathbf{I} \to \mathbf{R}$ also goes from $x$ to $x'$; since $a \leq c$, the affine interpolation from $a$ to $c$ is an extended 2-homotopy $a \prec_2 c$ (3.4); similarly, $b \prec_2 c$ and $a \simeq_2 b$, which projects to $[\bar{a}] = [\bar{b}]$. $\qquad\square$

**7.3. Theorem A** (Isometric classification). The w-spaces $w\mathbf{R}/G_\vartheta$ and $w\mathbf{R}/G_{\vartheta'}$ are isometrically isomorphic if and only if $G_\vartheta^+ = G_{\vartheta'}^+$ (as subsets of $\mathbf{R}$), if and only if $G_\vartheta = G_{\vartheta'}$, if and only if $\vartheta' \in \mathbf{Z} \pm \vartheta$.

**Proof.** If our w-spaces are isometrically isomorphic, also their fundamental weighted monoids (independently of the base point) are: $G_\vartheta^+ \cong G_{\vartheta'}^+$ (isometrically). Since the values of the weight $w: G_\vartheta^+ \to \mathbf{R}$ form the set $G_\vartheta^+$, it follows that $G_\vartheta^+ = G_{\vartheta'}^+$, which implies that $G_\vartheta$ (the additive subgroup of $\mathbf{R}$ generated by $G_\vartheta^+$) coincides with $G_{\vartheta'}$. If this is the case, then $\vartheta = a + b\vartheta'$ and $\vartheta' = c + d\vartheta$ for suitable integers  a, b, c, d;  whence  $\vartheta = a + bc + bd\vartheta$  and  $d = \pm 1$, so that  $\vartheta' = c \pm \vartheta$.  Finally, if $\vartheta' \in \mathbf{Z} \pm \vartheta$,  then  $G_\vartheta = G_{\vartheta'}$  and  $w\mathbf{R}/G_\vartheta = w\mathbf{R}/G_{\vartheta'}$. □

**7.4. Theorem B** (Lipschitz classification). The w-spaces  $w\mathbf{R}/G_\vartheta$  and  $w\mathbf{R}/G_{\vartheta'}$  are Lipschitz isomorphic if and only if the equivalent conditions of the following Lemma hold.

**Proof.** One implication follows from Theorem 7.2: if our w-spaces are Lipschitz isomorphic, also their fundamental weighted monoids  $G_\vartheta^+$  and  $G_{\vartheta'}^+$  are, by the functorial properties of  $w\Pi_1$  (3.4).

For the converse, let  $\vartheta'$  belong to the closure  $\{\vartheta\}_{RT}$; it suffices to consider the cases  $\vartheta' \in \vartheta + \mathbf{Z}$ and  $\vartheta' = \vartheta^{-1}$.  In the first case,  $G_\vartheta$  and  $G_{\vartheta'}$  coincide, as well as their action on  $w\mathbf{R}$; in the second, the Lipschitz isomorphism of weighted spaces

(1)   $f: w\mathbf{R} \to w\mathbf{R}$,                          $f(t) = |\vartheta|.t$,

restricts to an isomorphism  $f': G_\vartheta \to G_{\vartheta'}$,  obviously consistent with the actions  $(f(t + g) = f(t) + f'(g))$,  and induces a Lipschitz isomorphism  $w\mathbf{R}/G_\vartheta \to w\mathbf{R}/G_{\vartheta'}$. □

**7.6. Lemma.** Let  $\vartheta, \vartheta'$  be irrationals. The following conditions are equivalent:

(a)   the weighted groups  $G_\vartheta$  and  $G_{\vartheta'}$  are Lipschitz isomorphic,

(b)   the weighted monoids  $G_\vartheta^+$  and  $G_{\vartheta'}^+$  are Lipschitz isomorphic,

(c)   $G_\vartheta$  and  $G_{\vartheta'}$  are isomorphic as ordered groups (with respect to the total orders induced by  $\mathbf{R}$),

(d)   $\vartheta$  and  $\vartheta'$  are conjugate under the action of  $GL(2, \mathbf{Z})$  (7.1),

(e)   $\vartheta'$  belongs to the closure  $\{\vartheta\}_{RT}$  of  $\{\vartheta\}$  under the mappings  $R(t) = t^{-1}$  and  $T^{\pm 1}(t) = t \pm 1$.

**Proof.** The equivalence of the last three conditions is well-known, within the classification of the C*-algebras  $A_\vartheta$  up to strong Morita equivalence. It is also proved in [G5], Lemma 4.7.

Further, (a) implies (b), because  $G_\vartheta^+$  is the monoid of elements of  $G_\vartheta$  having a finite weight. And (b) implies (c), because  $G_\vartheta$  is the group canonically associated to the cancellative monoid  $G_\vartheta^+$, ordered with the latter as a positive cone.

Finally, to prove that (e) implies (a), let  $\vartheta'$  belong to the closure  $\{\vartheta\}_{RT}$; again, it suffices to consider the cases  $\vartheta' \in \vartheta + \mathbf{Z}$  and  $\vartheta' = \vartheta^{-1}$.  In the first,  $G_\vartheta = G_{\vartheta'}$; in the second, the Lipschitz isomorphism of weighted spaces  $f: w\mathbf{R} \to w\mathbf{R}$  considered above (7.4.1) restricts to a Lipschitz isomorphism of weighted abelian groups  $G_\vartheta \to G_{\vartheta'}$. □